\documentclass[11pt]{article}

\usepackage{amssymb,amsmath,amsthm,amsopn,amsfonts}
\usepackage[T1]{fontenc}
\usepackage{ae,aecompl}

\usepackage{graphicx}

\usepackage{bm}
\usepackage{fancybox,layout,color}
\usepackage{fancyhdr}

\usepackage{scrtime}

\usepackage{setspace}
\usepackage{url}
\usepackage{float}
\usepackage{placeins}

\numberwithin{equation}{section}
\definecolor{purple}{RGB}{160,32,40}

\usepackage{ulem}  

\newtheorem{teo}{Theorem}[section]

\newtheorem{defi}[teo]{Definition}
\newtheorem{prop}[teo]{Proposition}
\newtheorem{lema}[teo]{Lemma}
\newtheorem{coro}[teo]{Corollary}
\newtheorem{nota}[teo]{Remark}

\newtheorem{ex}[teo]{Example}

\newcommand{\Epi}{\mathrm{epi}}

\DeclareMathOperator{\co}{\mathsf{co}}
\DeclareMathOperator{\R}{\mathbb{R}}

\newcommand{\be}{\begin{eqnarray}}
\newcommand{\ee}{\end{eqnarray}}

\newcommand{\ben}{\begin{eqnarray*}}
\newcommand{\een}{\end{eqnarray*}}

\newcommand{\dist}{\mathrm{dist}}

\begin{document}

\title{Compensated Convex Transforms and\\  Geometric
Singularity Extraction from Semiconvex Functions\footnote{A Chinese version of the material in this manuscript has been published in 
{\it Zhang, Kewei, Crooks, Elaine and Orlando, Antonio, Compensated convex transforms and geometric singularity extraction from semiconvex functions (in Chinese), Sci. Sin. Math., 46 (2016) 747-768, doi: 10.1360/N012015-00339}}}

\author{\normalsize  Kewei Zhang\thanks{School  of Mathematical Sciences,
University of Nottingham, University Park, Nottingham, NG7 2RD, UK},\,
Elaine Crooks\thanks{Department of Mathematics, Swansea University,
Singleton Park, Swansea, SA2 8PP, UK}
\,
and Antonio Orlando\thanks{CONICET, Inst. de Estructuras \& Dept. de Mec\'{a}nica, 
Universidad Nacional de Tucum\'an, Argentina}
}

\date{{\normalsize \bf Dedicated to Professor Kung-ching Chang on the occasion of his 80th Birthday}}

\maketitle


\singlespacing
\pagestyle{fancy}
\fancyhead{}
\fancyfoot[OR,ER]{\tiny \today\,\, \thistime}
\cfoot{\thepage}


\begin{abstract}
	The upper and lower compensated convex transforms \cite{Z,Z2,ZOC} are `tight' one-sided approximations for a given function.  
	We apply these transforms to the extraction of fine geometric singularities from general semiconvex/semiconcave functions and
	DC-functions in $\mathbb{R}^n$ (difference of convex functions). Well-known geometric examples of (locally) semiconcave functions include the 
	Euclidean distance function and the Euclidean squared-distance function. 
	For a locally semiconvex function $f$ with general modulus, we show that `locally' a point is singular
	 (a non-differentiable point) if and only if it is a scale $1$-valley point, hence by using our method we can extract all fine singular points
	 from a given semiconvex function. More precisely, if $f$ is a semiconvex function with general modulus and $x$ is a singular point, then  
	 locally the limit of the scaled valley transform exists at every point $x$ and  can be calculated as
	$\displaystyle \lim_{\lambda\to+\infty}\lambda V_\lambda (f)(x)=r_x^2/4$, where $r_x$ is the radius of the minimal bounding sphere
	\cite{Jung} of the (Fr\'echet) subdiffential $\partial_- f(x)$ of the locally semiconvex $f$
	and $V_\lambda (f)(x)$ is the valley transform at $x$.
	Thus the limit function $\displaystyle \mathcal{V}_\infty(f)(x):=\lim_{\lambda\to+\infty}\lambda V_\lambda (f)(x)=r_x^2/4$ provides
	a `scale $1$-valley landscape function' of the singular set for a locally  semiconvex function $f$. 
	At the same time, the limit also provides an asymptotic expansion of the upper transform $C^u_\lambda(f)(x)$ when 
	$\lambda$ approaches $+\infty$.
	For a locally semiconvex function $f$ with linear modulus we show further that the limit of the gradient of the upper
	compensated convex transform $\displaystyle \lim_{\lambda\to+\infty}\nabla C^u_\lambda(f)(x)$ exists and equals the 
	centre of the minimal
	bounding sphere of $\partial_- f(x)$.
	We also show that for a DC-function $f=g-h$, the scale $1$-edge transform, when $\lambda\to+\infty$, satisfies
	$\displaystyle \liminf_{\lambda\to+\infty}\lambda E_\lambda (f)(x)\geq (r_{g,x}-r_{h,x})^2/4$, where $r_{g,x}$ and $r_{h,x}$ are the
	radii of the minimal bounding spheres of the subdifferentials $\partial_- g$ and $\partial_- h$ of the two convex functions
	$g$ and $h$ at $x$, respectively.
\end{abstract}

\medskip

\footnotesize
{\bf Keywords}:\textit{Compensated convex transforms, ridge transform, valley transform, edge transform, convex function,
	semiconvex function, semiconcave function, linear modulus, general modulus, DC-functions, singularity extraction,
	minimal bounding sphere,  local approximation, local regularity, singularity landscape}

\medskip
{\bf 2000 Mathematics Subjects Classification number}: 52A41, 41A30, 26B25, 49J52

\medskip
{\bf Email}:  kewei.zhang@nottingham.ac.uk,  e.c.m.crooks@swansea.ac.uk, aorlando@herrera.unt.edu.ar

\normalsize

\bigskip
\setcounter{equation}{0}
\section{Introduction and main results}

\bigskip
About ten years ago, the first author submitted the paper \cite{Z} dedicated to Professor Kung-ching Chang
on the occasion of his 70th birthday. Ten years on, the subject discussed in \cite{Z} 
has seen some further theoretical
developments \cite{Z2,ZOC,ZOCb,ZOCc}. As a step towards applications, we have been granted a UK patent \cite{ZOCd} 
on image processing methods based on this theory. 
In the present paper we work along a similar line to that in \cite{Z}.
We study the approximations and geometric singular extractions for semiconvex and semiconcave functions
by using compensated convex transforms introduced in \cite{Z}.

\medskip
Semiconcave and semiconvex functions have been extensively studied in the context of Hamilton-Jacobi equations \cite{CS}.
DC-functions (difference of convex functions) \cite{Ha} have been used in many optimisation problems \cite{UL}.
Important classes of such functions include the Euclidean distance function and the squared distance function.
Since general DC-functions and
semiconvex/semiconcave functions are locally Lipschitz functions in their essential domains (\cite[Theorem 2.1.7]{CS}),
 Rademacher's theorem implies that they are therein differentiable almost everywhere.
 Fine properties for the singular sets of convex/concave and semiconvex/semiconcave functions have been
 studied extensively \cite{AAC,AC,CS} showing that the singular set of a semiconvex/semiconcave
 function is rectifiable. However, from the applied mathematics point of view, natural questions arise,
 such as how such functions can be effectively approximated by smooth functions, whether all singular
 points are of the same type, that is, for semiconcave (semiconvex) functions, whether all singular
 points are geometric `ridge' (`valley') points, how singular sets can be effectively extracted beyond
 the definition of differentiability and how the information concerning `strengths' of different
 singular points can be effectively measured.
 Answers to these questions have important applications in image processing and computer-aided
 geometric design. For example,  the singular set of the Euclidean squared-distance
 function $\dist^2(\cdot,\Omega^c)$ to the complement of a bounded open domain
 $\Omega\subset\mathbb{R}^n$ (called the medial axis \cite{Bl} of the domain $\Omega$)
 carries important `compact' geometric information of the domain. It is also well known
 that the squared Euclidean distance function $\dist^2(\cdot,K)$ is $2$-semiconcave \cite{CS}.
 An answer to the question of how to extract the medial axis in a `stable' manner with respect
 to the domain under consideration has been addressed in \cite{ZOCc} and  has many applications \cite{SP}.
 In \cite{ZOCc} we introduced
 the notion of the {\it medial axis map} defined by
$M_\lambda(K)(x)=(1+\lambda)R_\lambda(\dist^2(\cdot, K))(x)$ for a closed set $K\subset\R^n$,
where  $R_\lambda(f)$ is the ridge transform of $f$ defined in \cite{ZOC}, and studied its properties.
We showed that $M_\lambda(\Omega)$ defines a Hausdorff stable multiscale representation of the
medial axis for finite $\lambda>0$ and the limit
$\lim_{\lambda\to+\infty}M_\lambda(\Omega)(x)=\dist^2(x,K)-\dist^2(x,\co[K(x)])$ exists for
all $x\in\R^n$, where $K(x)=\{y\in K,\; \dist(x,K)=|x-y|\}$ and $\co[K(x)]$ is its convex hull. 
This provides a `multiscale landscape' of the medial axis in the sense that higher
is the height, higher is the distance between the generating points of the medial axis branch.

\medskip
The present work is partly motivated  by \cite{ZOCc}. Our approximation results in the present
work are much more general than those in \cite{ZOCc}. Simple examples
which were not covered in \cite{ZOCc} are the Euclidean distance function itself and the weighted
squared distance function \cite{OBSC} for a finite set $K=\{x_i,\;i=1,\dots,m\}$ defined by
$\dist^2_{w,b}(x,\, K)=\min\{ w_i|x-x_i|^2+b_i,\; x_i\in K,\, w_i>0,\, b_i\in\R\}$. It is known
that the Euclidean distance function $\dist(\cdot,\,K)$  is locally semiconcave of linear modulus
in $\R^n\setminus K$ \cite{CS} and its singular set is more difficult to study geometrically than
that of the squared Euclidean distance function. It can be easily verified that the weighted squared
distance function is globally semiconcave. However, singularities for both of these functions are
difficult to study at a `finite scale'. This is in contrast with the standard Euclidean
functions \cite{ZOCc}.

\medskip
In \cite{ZOC,ZOCb}, we introduced several singularity extraction devices for detecting
geometric ridges, valleys, edges for functions and geometric intersections between smooth
manifolds defined by their characteristic functions (point clouds) based on compensated
convex transforms. These tools can also be used to measure the strength of singularities of
a particular type at a finite scale. In this paper we apply these tools to extract fine geometric
singularities from semiconvex/semiconcave functions and from DC-functions.
Our results demonstrate that our tight approximations by compensated convex transforms are of
very high quality in the sense they can extract geometric information of the original semiconvex/semiconcave
functions up to the first order derivative.

\medskip
We denote by $\R^n$ the standard $n$-dimensional Euclidean space with standard inner product
$x\cdot y$ and norm $|x|$ for $x,\,y\in\mathbb{R}^n$. We denote by $\bar A$  the closure of
a set $A$ in $\R^n$ and by  $B_r(x)$ and $\bar B_r(x)$ the open and closed balls in  $\R^n$ centred at
$x\in\R^n$ with radius $r>0$. We also denote by $C^1(\bar B_r(x))$ the space of real-valued continuously 
differentiable functions in an open set containing $\bar B_r(x)$ and by $C^{1,1}(\bar B_r(x))$
the space of real-valued continuous differentiable functions whose gradients are Lipschitz mappings.
Before we state our main results, let us first introduce the notions of compensated convex
transforms in $\mathbb{R}^n$. We state the definitions only for functions of linear growth
which will cover functions we deal with in this paper. For definitions under more general
growth conditions, see \cite{Z}. Let $f:\mathbb{R}^n\mapsto \mathbb{R}$ satisfy the linear
growth condition $|f(x)|\leq C|x|+C_1$ for some constants $C\geq 0$ and $C_1>0$ and for all $x\in \R^n$.

\medskip
The lower compensated compensated convex transform (lower transform for short) (see \cite{Z})
for $f$ is defined for $\lambda>0$ by
\begin{equation}
	C^l_\lambda(f)(x)=\co[f+\lambda|\cdot|^2](x)-\lambda|x|^2, \quad x\in\mathbb{R}^n,
\end{equation}
where $\co[g]$ is the convex envelope \cite{R,H-UL} of a function $g:\mathbb{R}^n\mapsto (-\infty,\, +\infty]$,
whereas the upper compensated compensated convex transform
(upper transform for short) (see \cite{Z}) for $f$ is defined for $\lambda> 0$ by
\begin{equation}
	C^u_\lambda(f)(x)=\lambda|x|^2-\co[\lambda|\cdot|^2-f](x),\quad x\in\mathbb{R}^n.
\end{equation}

\medskip
The two mixed compensated convex transforms are defined by $C^u_\tau(C^l_\lambda)(f)$ and
$C^l_\tau(C^u_\lambda)(f)$ when $\lambda,\,\tau>0$.

\medskip
It is known \cite{ZOC} that the lower and upper transforms are respectively the critical mixed Moreau
envelopes \cite{Mor65,Mor66,LL86,AA} and they can be viewed as morphological openings and closings \cite{ZOC}
respectively, in mathematical morphology terms \cite{Ser82,Jac92}.

\medskip
Since our main aim is to describe the behaviour of the ridge,  valley and edge transforms
for large $\lambda>0$, we introduce the following local versions of compensated convex transforms.
Due to the `locality property' for compensated convex transforms (see Proposition \ref{Prop.Loc} below), 
it will be obvious later that such definitions do not depend on the choices of domains involved.

\medskip
Let $\Omega\subset \R^n$ be an open set and let $f:\Omega\mapsto\R$ be a locally Lipschitz function,
which is thus bounded on every compact subset of $\Omega$. Assume $x\in\Omega$ and
let $G$ be a bounded open subset of $\Omega$  such that $x\in G\subset \bar{G}\subset \Omega$.
Let $L_G\geq 0$ be the Lipschitz constant of $f$ restricted to $\bar G$ denoted
by $f|_{\bar G}:\bar G\mapsto\R$. By Kirszbraun's theorem \cite{EG}, $f|_{\bar G}$ can be
extended to $\R^n$ as a Lipschitz continuous function $f_G:\R^n\mapsto \R$ with the same Lipschitz
constant $L_G$. Of course such an extension is not unique. However, due to the locality property
of compensated convex transforms, our results are independent of the Lipschitz extensions given
by Kirszbraun's theorem and the choices of $G$.

\medskip
Now we define the local lower compensated convex transform (local lower transform for short) and the
local upper compensated convex transforms (local upper transform for short) for a locally
Lipschitz function $f:\Omega\mapsto \R$ at $x\in \Omega$ with respect to $G$ respectively by
\begin{equation}
	C^l_{\lambda,G}(f)(x)=C^l_{\lambda}(f_G)(x)
	\quad\text{and}\quad
	C^u_{\lambda,G}(f)(x)=C^u_{\lambda}(f_G)(x),\quad x\in\R^n.
\end{equation}

\medskip
In \cite{ZOC} we introduced the notions of the ridge transform $R_\lambda(f)$,
the valley $V_\lambda(f)$ transform and the edge transform $E_\lambda(f)$, respectively, as
\begin{equation}
\begin{split}
	&R_\lambda(f)(x)=f(x)-C^l_\lambda(f)(x),\qquad V_\lambda(f)(x)=C^u_\lambda(f)(x)-f(x),\\[1.5ex]
	&E_\lambda(f)(x)=C^u_\lambda(f)(x)-C^l_\lambda(f)(x)=R_\lambda(f)(x)+V_\lambda(f)(x)
\end{split}
\end{equation}
for $x\in\R^n$.

\medskip
We should point out that our valley transform defined here is always non-negative 
and there is a sign difference in comparison with the valley transform defined in \cite{ZOC}.  
Given an open set $\Omega\subset \R^n$ and a locally Lipschitz function $f:\Omega\mapsto \R$,
we also define the local versions of ridge, valley and edge transforms as follows.

\medskip
\begin{defi}
	For $x\in\Omega$ and for a fixed open set $G$ whose closure is compact and
	$G$ satisfies $x\in G\subset \bar G\subset\Omega$, we define the local ridge, valley and edge
	transforms of $f$ at $x$ with respect to $G$ respectively as
\begin{equation}
	R_{\lambda,G}(f)(x)=R_{\lambda}(f_G)(x),\quad V_{\lambda,G}(f)(x)= V_{\lambda}(f_G)(x),
	\quad E_{\lambda,G}(f)(x)= E_{\lambda}(f_G)(x).
\end{equation}
\end{defi}

\medskip
Suppose $f:\mathbb{R}^n\mapsto \mathbb{R}$ is a Lipschitz function with Lipschitz constant $L\geq 0$.
It was established in \cite[Theorem 2.12 (iii)]{ZOC} that
\begin{equation}
	C^l_\lambda(f)(x)\leq f(x)\leq C^l_\lambda(f)(x)+\frac{L^2}{4\lambda},\quad
	C^u_\lambda(f)(x)-\frac{L^2}{4\lambda}\leq f(x)\leq C^u_\lambda(f)(x),
\end{equation}
for $\lambda>0$. Hence, the following estimates also hold \cite{ZOC}
\begin{equation}\label{Eq.EstLip}
	0\leq R_\lambda(f)(x)\leq \frac{L^2}{4\lambda},\qquad
	0\leq V_\lambda(f)(x)\leq \frac{L^2}{4\lambda}
 \end{equation}
 for $\lambda>0$, and at every point $x_0\in\R^n$  where $f$ is differentiable, we have
 \begin{equation}\label{Eq.AsyExt}
	\lim_{\lambda\to\infty}\lambda R_\lambda(f)(x_0)=0\quad\text{and}\quad
	\lim_{\lambda\to\infty}\lambda V_\lambda(f)(x_0)=0,
	\quad\text{hence}\quad
	\lim_{\lambda\to\infty}\lambda E_\lambda(f)(x_0)=0.
 \end{equation}
 For convenience later we call the quantities $\lambda R_\lambda(f)$, $\lambda V_\lambda(f)$ and $\lambda E_\lambda(f)$
 the scale $1$-ridge, -valley and -edge transforms, respectively.

\medskip
We will need also the following result
on the minimal bounding sphere for a compact set in $\R^n$.
The question  was first asked by J. J. Sylvester in  a two line statement \cite{Sy1} in 1857 for  finite sets in the plane,
which he then studied in his 1860 paper \cite{Sy2}.
The general result was proved by Jung in 1901 \cite{Jung}. There are however many later elementary proofs 
\cite{BW,V,DGK} by using Helly's theorem \cite{He}.

\medskip
\begin{lema}\label{Lem.BndSph} 
	(\cite{Jung, BW,V,DGK})
	Let $K\subset\R^n$ be a non-empty compact set. Then
	\begin{itemize}
		\item[$(i)$]	There is a unique minimal closed ball $\bar B_r(y_0)$ containing $K$ in the sense that
				$\bar B_r(y_0)$ is the closed ball  containing $K$ with the smallest radius.
				The sphere $S_r(x_0):=\partial B_r(x_0)$ is called the minimal bounding sphere of $K$.
		\item[$(ii)$]	Let $d$ be the diameter of $K$, then $r\leq \sqrt{\frac{n}{2(n+1)}}d$.
		\item[$(iii)$]	The centre of the ball $x_0$ satisfies $x_0\in \co[K\cap S_r(x_0)]$, the convex hull of $K\cap S_r(x_0)$.
	\end{itemize}
\end{lema}

\medskip
The proofs of Lemma \ref{Lem.BndSph}$(i)$ and $(ii)$ can be found in \cite{BW}
while for the proof of $(iii)$ we refer to \cite[2.6 and 6.1]{DGK}  or \cite[Lemma 2]{FG}.

\medskip
In this paper we will consider semiconvex and semiconcave functions, which are defined as follows \cite{CS,A}

\medskip
\begin{defi}\label{Def.SemCnvx}
	Let $\Omega\subset\mathbb{R}^n$ be  a non-empty open convex domain.
	\begin{itemize}
	\item[$(i)$] A function $f:\Omega\mapsto \mathbb{R}$ is called {\bf semiconvex} in $\Omega$ with modulus $\omega$ if there is a
		non-decreasing upper semicontinuous function $\omega:[0,\, +\infty)\mapsto [0,\,+\infty)$ such that
		$\displaystyle \lim_{t\to 0+}\omega(t)=0$
		and
		\begin{equation}
			sf(x)+(1-s)f(y)-f(s x+(1-s)y)\geq -s(1-s)|x-y|\omega(|x-y|)
		\end{equation}
		for all $x,\, y\in \Omega$ and for all $0\leq s\leq 1$.
	\item[$(ii)$] A function $f:\Omega\mapsto \mathbb{R}$ is {\bf semiconcave} in $\Omega$ with modulus $\omega$ if
		$-f$ is semiconvex with modulus $\omega$.
	\item[$(iii)$]	When $\omega(r)=\lambda_0r$ for $r\geq 0$ and for some $\lambda_0\geq 0$, we say that
			$f:\Omega\mapsto \mathbb{R}$ is {\bf $2\lambda_0$-semiconvex with linear modulus} \cite{CS} 
			($2\lambda_0$-semiconvex for short). In this case, there is a convex function 
			$g:\Omega\mapsto\R$ such that $f(x)=g(x)-\lambda_0|x|^2$ for all $x\in\Omega$ \cite[Propostion 1.1.3]{CS}.\\
 A function $f$ is {\bf $2\lambda_0$-semiconcave
			with linear modulus} ($2\lambda_0$-semiconcave for short)
			if $-f$ is $2\lambda_0$-semiconvex with linear modulus. 
			In this case, there is a concave function $g:\Omega\mapsto\R$ such that $f(x)=g(x)+\lambda_0|x|^2$ 
			for all $x\in\Omega$ \cite[Propostion 1.1.3]{CS}.
	\item[$(iv)$]	A function $f:\Omega\mapsto \mathbb{R}$ is called {\bf locally semiconvex} (respectively, {\bf locally semiconcave})
		in $\Omega$ if, on every convex compact set $K\subset\Omega$, $f$ is semiconvex (respectively, semiconcave)
		with a modulus $\omega_K$ depending on $K$.
\item[$(v)$] A function $f:\Omega\mapsto \mathbb{R}$ is called locally semiconvex (respectively, locally semiconcave) with linear modulus 
	if for every convex compact subset $K\subset \Omega$, there is a constant $\lambda_K\geq 0$ and a convex function ( respectively, concave function) 
	$g_K:K\mapsto\R$ such that when $x\in K$, we have $f(x)=g_K(x)-\lambda_K|x|^2$ (respectively, $f(x)=g_K(x)+\lambda_K|x|^2$).

	\end{itemize}
\end{defi}

\medskip
From Definition \ref{Def.SemCnvx}, it can be easily seen that the lower and upper compensated
convex transforms with scale $\lambda>0$ are $2\lambda$-semiconvex and $2\lambda$-semiconcave functions, respectively.
In fact, they are $2\lambda$-semiconvex and $2\lambda$-semiconcave `envelopes' of the given function.

\medskip
Let $\Omega\subset\R^n$ be a non-empty open convex set. We also recall  \cite[pag. 221]{A} that a locally semiconvex/semiconcave function $f:\Omega\to\R^n$ is locally
Lipschitz continuous in $\Omega$, that is, in every compact subset $K\subset \Omega$, $f$ is a Lipschitz function on $K$.

\medskip
The following is our main result on local approximations and geometric singular extraction of
semiconvex functions by the upper transform. The result regards
the Fr\'echet subdifferential of semiconvex functions. For its definition, we refer to Definition \ref{Def.SubDff} below
and to its characterization \eqref{Eq.DefSubSC}.

\medskip
\begin{teo}\label{Thm.LimV}  
	\begin{itemize}
		\item[$(i)$]	Let $\Omega \subset\R^n$ be a non-empty  open convex domain.
				Suppose $f:\Omega\mapsto \mathbb{R}$ is a locally semiconvex function in $\Omega$.
				Let $x_0\in\Omega$ be a non-differentiable (singular) point of $f$.
				Then for every bounded open set $G\subset \Omega$
				such that $x_0\in G\subset\bar{G}\subset \Omega$, 
				\begin{equation}
					\lim_{\lambda\to+\infty}\lambda V_{\lambda,G}(f)(x_0)=\frac{r_{x_0}^2}{4},
				\end{equation}
				where $r_{x_0}>0$ is the radius of the minimal bounding sphere of the subdifferential
				$\partial_- f(x_0)$ of $f$ at $x_0$.
		\item[$(ii)$]	\textcolor{black}{Assume that $f:\Omega\to \mathbb{R}$ is a locally semiconvex function with linear modulus in $\Omega$, i.e. on every convex compact subset
$K$ of $\Omega$, there exists $\lambda_K\geq 0$ such that $f(x)=g_K(x)-\lambda_K|x|^2$ for $x\in K$, where $g_K:K\to\mathbb{R}$ is a convex continuous function on $K$,
				 and  let $x_0\in\Omega$ be a non-differentiable (singular) point of $f$.}
				Then for every bounded open set $G\subset \Omega$ such that $x_0\in G\subset\bar{G}\subset \Omega$, 
				\begin{equation}
					\lim_{\lambda\to+\infty} \nabla C^u_{\lambda,G}(f)(x_0)=y_0,
				\end{equation}
				where $y_0\in \partial_- f(x_0)$ is the centre of the minimal bounding sphere of $\partial_-f(x_0)$.
	\end{itemize}
\end{teo}

\medskip
A similar result holds also for locally semiconcave functions, with the differences that we have to replace the
valley transform by the ridge transform so that $(i)$ of Theorem \ref{Thm.LimV} reads
\begin{equation}\label{Eq.SclRdg}
	\lim_{\lambda\to+\infty}\lambda R_{\lambda,G}(f)(x_0)=\frac{r_{x_0}^2}{4},
\end{equation}
with $r_{x_0}>0$ the radius of the minimal bounding sphere of the (Fr\'echet) superdifferential $\partial_+ f(x_0)$
of the locally semiconcave function $f$ at $x_0$ (see Definition \ref{Def.SupDff} below), while $(ii)$ becomes
\begin{equation}\label{Eq.GrdLwTr}
	\lim_{\lambda\to+\infty} \nabla C^l_{\lambda,G}(f)(x_0)=y_0,
\end{equation}
with $y_0\in \partial_+ f$ the centre of the minimal bounding sphere of $\partial_+f(x_0)$.

\medskip
Since near every point
$x\in G$, with $G$ a bounded open subset of $\Omega$ such that $x\in G\subset \bar{G}\subset \Omega$,
$C^u_\lambda(f_G)$ is a $C^1$ function in any given neighbourhood
$B_r(x)\subset\bar B_r(x)\subset G$ for sufficiently large $\lambda>0$ due to the locality property (see Proposition \ref{Prop.Loc} below),
$C^u_\lambda(f_G)$ realizes a locally smooth approximation from above and the error of the
approximation satisfies
\[
	\lambda V_\lambda(f_G)(x)=\lambda(C^u_\lambda(f_G)(x)-f_G(x))\to r_x^2/4\quad\text{for }\lambda\to +\infty
\]
at a singular point $x\in G$.

\medskip
In order to help readers to have an intuitive view on compensated convex transforms,
the ridge/valley transforms and their limit for semiconvex/semiconcave functions, we consider the
following simple  example first.

\medskip
\begin{ex}
	Let $f(x)=|x|$ for $x\in\R$. Clearly, $f$ is a convex function. For $\lambda>0$, we have
\begin{equation}\label{Eq.1dEX}
	\begin{array}{ll}
			\begin{array}{l}
	\displaystyle C^u_\lambda(f)(x)=
		\left\{\begin{array}{l}
			\displaystyle	\lambda x^2+\frac{1}{4\lambda}, \quad |x|\leq \frac{1}{2\lambda},\\
			\displaystyle	|x|,\qquad |x|\geq \frac{1}{2\lambda};
		\end{array},\right.\\[2ex]
	\displaystyle \lim_{\lambda\to+\infty}\frac{d}{dx}C^u_\lambda(f)(x)=
		\left\{\begin{array}{ll}
				\displaystyle	-1,& x<0,\\
				\displaystyle	0,& x=0,\\
				\displaystyle	1,& x>0.
			\end{array}\right.
			\end{array}
			&
			\begin{array}{l}
		\displaystyle \lambda V_\lambda(f)(x)=
			\left\{\begin{array}{ll}
				\displaystyle	\lambda^2\left(|x|-\frac{1}{2\lambda}\right)^2,	& |x|\leq \frac{1}{2\lambda},\\
				\displaystyle	0,						& |x|\geq \frac{1}{2\lambda};
				\end{array}\right.\\[2ex]
	\displaystyle	\lim_{\lambda\to+\infty}\lambda V_\lambda(f)(x)=
			\left\{\begin{array}{ll}
				\displaystyle \frac{1}{4},& x=0,\\
				\displaystyle	0,& x\neq 0;
				\end{array}\right.
			\end{array}
	\end{array}
\end{equation}
For this example the subdifferential of $f$ at $0$ is given by $\partial_-f(0)=[-1,\,1]$.
Thus the smallest closed interval which contains $\partial_-f(0)$ coincides with $\partial_-f(0)$ itself, 
with the mid point $0$ and radius $1$.
Note also that Theorem \ref{Thm.LimV}$(i)$ and $(ii)$ hold
in this case.\hfill\qed
\end{ex}

\medskip
There are many examples of locally semiconvex/semicocave functions \cite{CS}.
Suppose $\Omega\subset \R^n$ is open and $K\subset \R^m$ is compact. If $F:K\times\Omega\mapsto \R$
and $\nabla_x F$ are both continuous in $K\times\Omega$, then $f(x)=\sup_{s\in K}F(s,x)$ is locally semiconvex.
If $\nabla^2_xF$ also exists and is continuous in $K\times\Omega$, then $f$ is locally semiconvex with
linear modulus (see \cite[Proposition 3.4.1]{CS}).

\medskip
The following are two important  examples on extraction of geometric singular points arising from applications.
They refer to the square distance function and to the distance function to a closed set $K\subset \R^n$.

\medskip
\begin{ex}
Let $K\subset \R^n$ be a non-empty closed set, satisfying $K\neq \R^n$ and denote by
$\dist^2(\cdot,\,K)$ the squared Euclidean distance function to $K$.
Let $M_K:=\{y\in \R^n,\,\exists z_1,\,z_2\in K,\, z_1\neq z_2\; \dist(x,\, K)=|y-z_1|=|y-z_2|\}$
be the medial axis of $K$. It is known that $M_K$ is the singular set of $\dist^2(\cdot,\,K)$.
In \cite{ZOCc} we have the following Luzin type theorem. Let $\lambda>0$. If we define
\[
	V_{\lambda,K}:=\{x\in \R^n,\quad \lambda\dist(x,\,M_K)\leq \dist(x,\,K)\},
\]
then
\[
	\dist^2(x,\, K)=C^l_\lambda(\dist(\cdot,\,K))(x)
\]
for $x\in\R^n\setminus V_{\lambda,K}$ and
\[
	\bar M_K=\cap^\infty_{\lambda>0}V_{\lambda,K}\,.
\]	
As a result, we have \cite{ZOCc}
\[
	\dist^2(\cdot, \, K)\in C^{1,1}(\R^n\setminus V_{\lambda,K})
\]
and
\[
	|\nabla \dist^2(x,\,K)-\nabla \dist^2(y,\,K)|\leq 2\max\{1,\,\lambda\}|x-y|,\quad x,\; y\in \R^n\setminus V_{\lambda, K}\,.
\]
Since the proof of this result relies on the special geometric features of the squared Euclidean
distance function, in \cite{ZOCc} we have not been able to extend this result to more general semiconcave functions.
We have therefore defined the (quadratic) medial axis map as $M_\lambda(x,\, K)= (1+\lambda)R_\lambda(\dist^2(\cdot,\,K)(x)$
and proved that
\begin{equation}\label{Eq.LmMMA}
	\lim_{\lambda\to+\infty}M_\lambda(x,\, K)=\dist^2(x,\,K)-\dist^2(x,\,\co[K(x)]),
\end{equation}
where $\co[K(x)]$ is the convex hull of the compact set $K(x)=\{y\in K,\; \dist(x,K)=|x-y|\}$.

\medskip
We can now interpret the limit \eqref{Eq.LmMMA} by applying Theorem \ref{Thm.LimV}$(i)$.
Since $\displaystyle \lim_{\lambda\to+\infty}R_\lambda(\dist^2(\cdot,\,K)(x)=0$,
by the definition of $M_\lambda(x,\, K)$ we have
\[
	\lim_{\lambda\to+\infty}M_\lambda(x,\, K)=\lim_{\lambda\to+\infty}\lambda R_\lambda(\dist^2(\cdot,\,K)(x),
\]
where $\lambda R_\lambda(\dist^2(\cdot,\,K)(x)$ is our scale $1$-ridge transform. Now, for
$x\in M_K$, the superdifferential of $\dist^2(\cdot,\,K)$ at $x$ is given by
$\partial_+\dist^2(x,\,K)=\co\{ 2(x-y),\, y\in K(x)\}$ so that the square $r_x^2$ of the radius
of the minimum bounding sphere of $\partial_+\dist^2(x,\,K)$ is $4\Big(\dist^2(x,\,K)-\dist^2(x,\, \co[K(x)])\Big)$.
Thus $r_x^2/4$, which is the limit of the scale $1$-ridge transform (see \eqref{Eq.SclRdg}), 
is the same as $\dist^2(x,\,K)-\dist^2(x,\,\co[K(x)])$ (see \eqref{Eq.LmMMA}, \cite[Theorem 3.23]{ZOCc}).
\end{ex}

\begin{ex}
In this example, we consider the case of the Euclidean distance function $\dist(x,\, K)$ itself.
It is then known \cite[Proposition 2.2.2]{CS} that $\dist(\cdot,K)$ is locally semiconcave with
linear modulus in $\R^n\setminus K$.
Therefore if we consider the limit of the scale $1$-ridge transform, by Theorem \ref{Thm.LimV}
applied to semiconcave functions, we have (see \eqref{Eq.SclRdg})
\[
	\lim_{\lambda\to+\infty}\lambda R_\lambda(\dist(\cdot,K))(x)=r_x^2/4\;\text{and}\;
	\lim_{\lambda\to+\infty}\nabla C^l_\lambda(\dist(\cdot,K))(x)=y_x,\quad x\notin K,
\]
where $r_x$ is the radius of the minimal bounding sphere of the superdifferential
$\partial_+\dist(x,K)$ and $y_x$ is the centre of the
minimal bounding sphere. Since $\partial_+\dist(x,K)=\co\{(x-y)/|x-y|,\; \dist(x,K)=|x-y|\}$,
if we let $p_x\in \co[K(x)]$ be the unique closest point from $x$ to $\co[K(x)]$,  then we have
\begin{equation}\label{Eq.Rd}
	r^2_x=\frac{\dist^2(x,\,K)-\dist^2(x,\,\co[K(x)]}{\dist^2(x,\,K)},\quad y_x=\frac{p_x}{\dist(x,\, K)}\,.
\end{equation}
By comparing \eqref{Eq.LmMMA} and \eqref{Eq.Rd} we find that for $x\notin K$
\begin{equation*}
	\lim_{\lambda\to+\infty}\lambda R_\lambda(\dist(\cdot,K))(x)
			=\frac{\displaystyle \lim_{\lambda\to+\infty}\lambda R_\lambda(\dist^2(\cdot,K))(x)}{4\dist^2(x,\,K)},\\[1.5ex]
\end{equation*}
and		
\begin{equation*}
		\lim_{\lambda\to+\infty} \nabla C^l_\lambda(\dist(\cdot,K))(x)
			=\frac{\displaystyle \lim_{\lambda\to+\infty} \nabla C^l_\lambda(\dist^2(\cdot,K))(x)}{2\dist(x,\,K)}
\end{equation*}
whereas for $x\in K$, we have that $R_\lambda(\dist(\cdot,K))(x)=R_\lambda(\dist^2(\cdot,K))(x)=0$
as points in $K$ are minimum points of both the distance function and the squared distance function \cite{Z}.
We can conclude therefore that Theorem \ref{Thm.LimV} links the asymptotic behaviours of $C^l_\lambda(\dist(\cdot,K))(x)$ and
$C^l_\lambda(\dist^2(\cdot,K))(x)$, with the latter which is much easier to analyse \cite{ZOCc}. \hfill \qed
\end{ex}

\medskip
For DC-functions, that is, functions that can be represented as difference between two convex functions, 
we have the following sufficient condition for extracting edges.

\begin{coro}\label{Coro.EdgDC}
Let $\Omega\subset\R^n$ be a non-empty open convex set. Assume $g,\, h: \Omega\mapsto \R$ are finite continuous
convex functions in $\Omega$ and let $f(x)=g(x)-h(x)$ for $x\in \Omega$.
Take $x_0\in \Omega$ and $G\subset\Omega$ an open bounded set such that $x_0\in G\subset\bar{G}\subset \Omega$.
Let $r_{g,x_0}$ and $r_{h,x_0}$ be the radii of the minimal bounding spheres of $\partial_-g(x_0)$
and $\partial_-h(x_0)$, respectively. Then,
\begin{equation}\label{Eq.DCEdg}
	\liminf_{\lambda\to+\infty}\lambda E_\lambda f_G(x_0)\geq \frac{(r_{g,x_0}-r_{h,x_0})^2}{4}.
\end{equation}
\end{coro}

\medskip
\begin{nota}
It is easy to see that the lower bound in \eqref{Eq.DCEdg} is sharp. If we set  $g(x)=h(x)=|x|$ for $x\in \R$,
$f\equiv 0$, thus $r_{g,0}=r_{h,0}=1$ while $E_\lambda(f)(0)=0$ for all $\lambda>0$.
However, when $r_{g,x_0}=r_{h,x_0}$, there are simple examples that show that the left hand side of \eqref{Eq.DCEdg} may be strictly positive.
For example, if we let $F(x,y)=|x|-|y|$ in $\R^2$ and let $f(x)=|x|$, it is easy to see that
$E_\lambda(F)(x,y)= V_\lambda(f)(x)+V_\lambda(f)(y)$, hence by \eqref{Eq.1dEX}, we have
$\displaystyle \lim_{\lambda\to+\infty} \lambda E_\lambda(F)(0,0)=1/2>0$.
Note that if we write $f_1(x,y)=f(x)$ and $f_2(x,y)=f(y)$, we have $\partial_-f_1(0,0)=[-1,1]\times\{0\}$ while
$\partial_-f_2(0,0)=\{0\}\times[-1,1]$.
The minimal bounding sphere for both $\partial_-f_1(0,0)$ and $\partial_-f_2(0,0)$ is the unit sphere in $\R^2$,
thus $r_{f_1,0}=r_{f_2,0}$. In general, it would be rather technical to analyse the left-hand side of \eqref{Eq.DCEdg}
based on the subdifferentials $\partial_-g(x_0)$ and $\partial_-h(x_0)$ \cite{UL}.
We will not consider this case here.
\end{nota}

\medskip
We say that compensated convex transforms are `tight approximations' for a given function. Roughly speaking
for functions that are locally of class $C^{1,1}$ near $x_0$, then there is a finite $\Lambda>0$, such that 
$C^u_\lambda(f)(x_0)=f(x_0)=C^l_\lambda(f)(x_0)$ whenever $\lambda \geq \Lambda$ \cite[Theorem 2.3(iv)]{Z}. 
This implies that at a smooth point, the graph of the upper/lower transform is tightly attached to that of the original function from above/below.
If $f:\Omega\mapsto \R$ is locally a semiconvex/semiconcave function with linear modulus, 
where $\Omega\subset \R^n$ is a non-empty convex open set, then according to the well-known Alexandrov's theorem \cite{Ev,CS}, 
$f$ is twice differentiable almost everywhere in $\Omega$,  that is,
for almost every $x_0\in \Omega$, there is some $p\in\R^n$ and an $n\times n$ symmetric matrix $B$ such that
\begin{equation}\label{Eq.AlexTheo}
	\lim_{x\to x_0}\frac{f(x)-f(x_0)-p\cdot (x-x_0)-(x-x_0)\cdot B(x-x_0)}{|x-x_0|^2}=0\,.
\end{equation}

\medskip
We say that $x_0\in\Omega$ is an Alexandrov point if \eqref{Eq.AlexTheo} holds.

\medskip
\begin{prop}\label{Prop.TghAprxScnvx}
	Let $\Omega\subset\R^n$ be a non-empty open convex set.
	Suppose $f:\Omega\mapsto \mathbb{R}$ is a locally semiconvex/semiconcave function of linear modulus.  
	Assume $x_0\in\Omega$ and  $G$ a bounded open subset of $\Omega$ such that
	$x_0\in G\subset \bar G\subset \Omega$. If $x_0$ is an Alexandrov point,
	there is a constant $\Lambda>0$, such that when $\lambda\geq \Lambda$, we have
	\begin{equation}\label{Eq.TghtAprx}
		f(x_0)=C^u_\lambda(f_G)(x_0)=C^l_\lambda(f_G)(x_0),
	\end{equation}
	and
	\begin{equation}\label{Eq.TghtAprx.Drv}
		\nabla f(x_0)=\nabla C^u_\lambda(f_G)(x_0)=\nabla C^l_\lambda(f_G)(x_0).
	\end{equation}
\end{prop}

\medskip
\begin{nota}
\begin{itemize}
\item[(i)] For a locally semiconvex function $f$ with linear modulus, it is not difficult to show
	that by the locality property, for every fixed $x\in G$, when $\lambda>0$ is sufficiently large,
	$f(x)=C^l_\lambda(f_G)(x)$. The slightly more involved part
	is to show that also the upper transform $C^u_\lambda(f_G)(x)$
	attains the value $f(x)$ for a finite $\lambda>0$ at an Alexandrov point.
\item[(ii)] Theorem \ref{Thm.LimV}, Proposition \ref{Prop.TghAprxScnvx} and \eqref{Eq.AsyExt}
	provide a clearer picture on how compensated convex transforms approach a locally semiconvex
	function with linear modulus.
\item[(iii)] Since at every point $x\in G\subset\bar G\subset\Omega$,
	$\displaystyle  \lim_{\lambda\to+\infty}\lambda V_{\lambda,G}(f)(x)$
	and $\displaystyle \lim_{\lambda\to+\infty}\lambda R_{\lambda,G}(f)(x)$ exist, we can define the
	`valley landscape  map' and
	the `ridge landscape  map' for locally semiconcovex and
	locally semiconcave functions with general modulus, respectively, by
	\begin{equation} \label{Eq.Lim}
		\mathcal{V}_\infty(f)(x)=\lim_{\lambda\to+\infty}\lambda V_{\lambda,G}(f)(x),
		\quad \mathcal{R}_\infty(f)(x)=\lim_{\lambda\to+\infty}\lambda R_{\lambda,G}(f)(x),
	\end{equation}
	Due to the locality property, the limits \eqref{Eq.Lim} are independent of the choice of $G$.
\item[(iv)] From the definition of the `valley landscape  map' of a semiconvex function $f$,
	we can identify at least three distinct features:
	\begin{itemize}
		\item[(a)] $\lambda V_{\lambda,G}(f)(x)=0$ in finite time $\lambda>0$ if $x$ is an Alexandrov point;
		\item[(b)] If  $f$ is differentiable at $x$ and $\lambda V_{\lambda,G}(f)(x)> 0$ for all $\lambda>0$,
				then $\displaystyle  \lim_{\lambda\to+\infty}\lambda V_{\lambda,G}(f)(x)=0$;
		\item[(c)] If $f$ is not differentiable at $x$, then
				$\displaystyle  \lim_{\lambda\to+\infty}\lambda V_{\lambda,G}(f)(x)=r_x^2/4>0$.
	\end{itemize}	
	Therefore, for large $\lambda>0$, subject to the boundary effect for points near $\partial G$,
	the set $\{x\in G,\; \lambda V_{\lambda,G}(f)(x)>\epsilon\}$ for a fixed $\epsilon>0$
	contains both singular points of $f$ in $G$ and points of high curvature, that is,
	either $\nabla^2 f(x)$ does not exist or the largest eigenvalue of $\nabla^2 f(x)$ is very large.
\end{itemize}
\end{nota}

\medskip
In Section \ref{Sec.Prel}, we introduce some further preliminary results which are needed for the proofs of our main
results Theorem \ref{Thm.LimV} and Corollary \ref{Coro.EdgDC}. We prove our results in Section \ref{Sec.Prfs}.


\bigskip
\setcounter{equation}{0}
\section{Some preliminary results}\label{Sec.Prel}
In this section, we collect some basic properties of compensated convex transforms which will be needed in the following,
and refer to \cite{Z,ZOC,ZOCb} for proofs and details.

The ordering property of compensated convex transforms holds for $x\in\mathbb{R}^n$ and reads as
\[
	C^l_\lambda(f)(x)\leq C^l_\tau(f)(x)\leq f(x)\leq C^u_\tau(f)(x)\leq C^u_\lambda(f)(x),\quad \tau\geq \lambda.
\]
The upper and lower transform for functions $f:\R^n\mapsto\R$ with quadratic growth, i.e.
$|f(x)|\leq C(1+|x|^2)$ for $x\in \R^n$ and for a constant $C\geq 0$, are related to each other
when $\lambda>0$ is large enough by the following relation
\[
	C^l_\lambda(f)(x)=-C^u_\lambda(-f)(x)\,.
\]
If $f$ is a continuous function with quadratic growth, 
\[
	\lim_{\lambda\to\infty}C^l_\lambda(f)(x)=f(x),\quad
	\lim_{\lambda\to\infty}C^u_\lambda(f)(x)=f(x),\quad x\in\R^n.
\]
If $f$ and $g$ are both Lipschitz functions, then for $\lambda>0$ and $\tau>0$, we have
\begin{equation}\label{Eq.InLip}
	C^l_{\lambda+\tau}(f+g)\geq C^l_{\lambda}(f)+C^l_{\tau}(g),\quad
	C^u_{\lambda+\tau}(f+g)\leq C^u_{\lambda}(f)+C^u_{\tau}(g).
\end{equation}

\medskip
We recall from \cite{BKK} the following definition.
\begin{defi}
We say that $f:\R^n\mapsto \R$ is upper semi-differentiable at $x_0\in\mathbb{R}^n$ if there is
some $u\in\mathbb{R}^n$ such that
\[
	\limsup_{y\to 0}\frac{f(x_0+y)-f(x_0)-u\cdot y}{|y|}\leq 0\,.
\]
\end{defi}

\medskip
The following differentiability property \cite[pag 726]{KK}  and more generally \cite[Corollary 2.5]{BKK} is useful in
the proofs of our results.

\medskip
\begin{lema}\label{Lem.DfUp}
	Suppose $g:B_r(x_0)\mapsto \R$ is convex and $f:  B_r(x_0)\mapsto \R$
	is upper semi-differentiable at $x_0$, such that $g\leq f$ on $B_r(x_0)$ and $g(x_0)=f(x_0)$.
	Then $f$ and $g$ are both differentiable at $x_0$ and $\nabla f(x_0)=\nabla g(x_0)$.
\end{lema}

\medskip
Note that concave functions are upper semi-differentiable.

\medskip
We recall the following locality property of the compensated convex transforms for Lipschitz continuous functions.
A similar result for bounded functions was established in \cite{ZOC}.

\medskip
\begin{prop}\label{Prop.Loc} 
	Suppose $f:\R^n\mapsto \R$ is Lipschitz continuous with Lipschitz constant $L>0$.
	Let $\lambda>0$ and $x\in \R^n$.
	Then there exist $(\tau_i,y_i)\in \R\times \R^n$, $i=1, \ldots, n+1$, such that
\begin{equation}\label{Eq.LocCmp}
	\begin{split}
		\co[f+\lambda|(\cdot) -x|^2](x)&=\co_{\bar B_{r_\lambda}(x)}[f+\lambda|(\cdot)-x|^2](x)\\[1.5ex]
						 &:=
						\inf\left\{\sum^{n+1}_{i=1}\tau_i[f(y_i)+\lambda|y_i-x|^2]:
		 y_i\in \R^n,\,\tau_i\geq 0,\,
		 |y_i-x|\leq r_\lambda,\right.\\
		 &\phantom{xxxxxxxxxxxxxxxxxxxxxxxxxx}\left.\sum^{n+1}_{i=1}\tau_i=1,\, \sum^{n+1}_{i=1}\tau_iy_i=x\right\}
	\end{split}
\end{equation}
where $r_\lambda=(2+\sqrt{2})L/\lambda$.\\
Furthermore, there is an affine function $y\mapsto \ell(y)=a\cdot (y-x)+b$
for $y\in\R^n$ with $a\in\mathbb{R}^n$ and $b\in\R$ such that
\begin{itemize}
	\item[(i)]	$\ell(y)\leq f(y)+\lambda|y-x|^2$ for all $y\in \mathbb{R}^n$;
	\item[(ii)]	$\ell(x_i)= f(x_i)+\lambda|x_i-x|^2$ for $i=1,\ldots,n+1$;
	\item[(iii)]	$b=\ell(x)=\co[f+\lambda|(\cdot)-x|^2](x)$.
\end{itemize}
\end{prop}

\medskip
We call $\co_{\bar B_{r_\lambda}(x)}[\lambda|(\cdot)-x|^2+f](x)$ defined in \eqref{Eq.LocCmp} the local convex envelope
of $y\in\R^n\mapsto \lambda|y-x|^2+f(y)$ at $x$ in $\bar B_{r_\lambda}(x)$.

\medskip
\begin{nota}
\begin{itemize}
	\item[(i)] The locality property given in Proposition \ref{Prop.Loc} also applies to the compensated convex transforms.
	Due to the translation invariance property \cite{ZOC}, for every fixed $x_0\in\R^n$, we have
	\begin{equation}
	\begin{split}
		C^l_\lambda(f)(x)&=\co[f+\lambda|(\cdot) -x_0|^2](x)-\lambda|x-x_0|^2,\\
		C^u_\lambda(f)(x)&=\lambda|x-x_0|^2-\co[\lambda|(\cdot) -x_0|^2-f](x)\,,
	\end{split}
	\end{equation}
	thus, if we take $x_0=x$, we obtain
	\begin{equation}
		C^l_\lambda(f)(x)=\co[f+\lambda|(\cdot) -x|^2](x),\quad
		C^u_\lambda(f)(x)=-\co[\lambda|(\cdot) -x|^2-f](x)\,,
	\end{equation}
	and \eqref{Eq.LocCmp} can be used.
	\item[(ii)] A consequence of  \cite[Remark 2.1]{Z} is that if $f$ is continuous and with linear growth, then the infimun
	in the definition of the convex envelope of the function $y\in\R^n\mapsto \lambda|y-x|^2+f(y)$ at $y=x$
	is attained by some $\lambda_i>0$, $x_i\in\R^n$, $i=1,\dots,k$ with $2\leq k\leq n+1$ (see \cite{H-UL,R}), that is, 
	\[
		\co_{\bar B_{r_\lambda}(x)}[f+\lambda|(\cdot)-x|^2](x)=\sum^{k}_{i=1}\lambda_i[f(x_i)+\lambda|x_i-x|^2]
	\]
	with $|x_i-x|<r_{\lambda}$, $i=1,\dots,k$ with $2\leq k\leq n+1$, and $\sum^{k}_{i=1}\lambda_i=1,\, \sum^{n+1}_{i=1}\lambda_ix_i=x$.
\end{itemize}
\end{nota}

\medskip
The following lemma can be considered a special case of Theorem \ref{Thm.LimV}.

\medskip
\begin{lema}\label{Lem.Ab}
	Let $S\subset \R^n$ be a non-empty compact convex set, containing more that one element, and denote by
	$S_r(-a)$ the minimal bounding sphere of $S$ with radius $r>0$ and centre $-a\in \R^n$.
	Consider the sublinear function $\sigma:x\in\R^n\to\sigma(x)=\max\{p\cdot x,\; p\in S\}$.
	Then for a fixed $0\leq\epsilon<\min\{1,\, r\}$ and for $\lambda>0$, we have
	\begin{equation}\label{Eq.Lem2.4.1}
		C^u_\lambda(\sigma-\epsilon|\cdot|)(0)=\frac{(r-\epsilon)^2}{4\lambda},
	\end{equation}
	\begin{equation}\label{Eq.Lem2.4.2}
		\nabla C^u_\lambda(\sigma)(0)=-a;
	\end{equation}
	and for a fixed $0<\epsilon<\min\{1,\, r\}$
	\begin{equation}\label{Eq.Lem2.4.3}
		C^u_\lambda(\sigma+\epsilon|\cdot|)(0)\leq C^u_{(1-\epsilon)\lambda}(\sigma)(0)+C^u_{\epsilon\lambda}(\epsilon|\cdot|)(0)=
		\frac{r^2}{4(1-\epsilon)\lambda}+\frac{\epsilon}{4\lambda},
	\end{equation}
	where
	\begin{equation}\label{Eq.Lem2.4.4}
		C^u_{\epsilon\lambda}(\epsilon|\cdot|)(x)=
			\left\{\begin{array}{ll}
				\displaystyle \epsilon\lambda|x|^2+\frac{\epsilon}{4\lambda},	&\displaystyle |x|\leq \frac{1}{2\lambda},\\[2ex]
				\displaystyle \epsilon|x|,					&\displaystyle |x|\geq \frac{1}{2\lambda}.
			\end{array}\right.
	\end{equation}
\end{lema}

\medskip
We have also the following local $C^{1,1}$ result for the upper transform of locally
semiconvex functions with linear modulus.

\medskip
\begin{prop}\label{Prop.UpTrC11}
	Let $f:\mathbb{R}^n\mapsto\R$ be a Lipschitz continuous function with Lipschitz constant $L\geq 0$.
	Assume that for some $r>0$,  $f$ is $2\lambda_0$-semiconvex in the closed ball $B_{2r}(0)$,
	that is,  $f(x)=g(x)-\lambda_0|x|^2$ for $x\in\bar B_{2r}(0)$, where  $\lambda_0\geq 0$
	is a constant and $g:\bar B_{2r}(0)\mapsto\R$ is convex. Then for $\lambda\geq \lambda_0$
	sufficiently large, $C^u_\lambda(f)\in C^{1,1}(\bar B_r(0))$ and
	\begin{equation} \label{Eq.EstDer}
		|\nabla C^u_\lambda(f)(x)-\nabla C^u_\lambda(f)(y)|\leq 2\lambda|x-y|,\quad x,\,y\in \bar B_r(0).
	\end{equation}
\end{prop}

\medskip
\begin{nota}
	From the proof of Proposition \ref{Prop.UpTrC11} (and \cite[Theorem 4.1]{Z} with a Lipschitz constant less as sharp)
	we can derive that if $f:\mathbb{R}^n\mapsto\R$ is both Lipschitz continuous and convex, for example
	if $f(x)=\sigma(x)$ is the sublinear function \cite{H-UL} defined by
	\[
		\sigma:x\in\R^n\to\sigma(x)=\max\{x\cdot p,\; p\in S\}\,,
	\]
	where $S$ is compact and convex, the estimate \eqref{Eq.EstDer} holds globally in $\R^n$ with $\lambda_0=0$.
\end{nota}

\medskip
We conclude this section by recalling the definition and some properties of the subdifferential of convex and semiconvex
functions we need in our proofs.

\medskip
\begin{defi}
Let $\Omega\subset\R^n$ be a non-empty open convex set. Assume $f:\Omega\mapsto\R$ is convex and let $x\in\Omega$.
The subdifferential of $f$ at $x$, denoted by $\partial_-f(x)$, is the set of $u\in \R^n$
satisfying \cite{H-UL}
\[
	f(y)-f(x) - u\cdot (y-x)\geq 0\,, \quad \text{for all }y\in\Omega.
\]
\end{defi}
\medskip

The subdifferential $\partial_-f(x)$ is a non-empty, compact and convex subset of $\R^n$.
If we define the sublinear function \cite[Chapter D]{H-UL} $y\in\R^n\to\sigma_x(y):=\max\{ u\cdot y,\; u\in\partial_-f(x)\}$
then
\begin{equation}\label{Eq.DirDer}
	\lim_{h\to 0}\frac{f(x+h)-f(x)-\sigma_x(h)}{|h|}=0,
\end{equation}
where $\sigma_x(h)$ defines the directional derivative of $f$ at $x$ along $h\in\R^n$.

\medskip

Just like the convex case, locally semiconvex functions have a natural notion of generalized gradient given by the subdifferential.
This is defined as follows.
\medskip

\begin{defi}\label{Def.SubDff}
	Let $f:\Omega\mapsto \R^n$ be a locally semiconvex function in $\Omega$ and let $x\in \Omega$. Denote by $K$
	an open convex subset of $\Omega$ such that $x\in K\subset \bar{K}\subset \Omega$ and by $\omega_K$
	a semiconvex modulus for $f$ in $K$.
	The Fr\'echet subdifferential $\partial_-f$ of $f$ at $x$ is the set of vectors $p\in \R^n$ satisfying
	\begin{equation}\label{Eq.DefSubSC}
		f(y)-f(x)- p\cdot (y-x)	\geq -|y-x|\omega_K(|y-x|)
	\end{equation}
	for any point $y$ such that the segment of ends $y$ and $x$ is contained in $K$.
\end{defi}

\medskip
It is not difficult to show that the definition	of $\partial_-f(x_0)$ does not depend on $K$, in fact,
condition \eqref{Eq.DefSubSC} can be expressed in terms of a kind of regularization of the semiconvexity modulus (see \cite[Proposition 2.1]{A}).
We also have that $\partial_-f(x_0)$ is a non-empty convex compact set.
Likewise for convex functions, we can equally define for locally semiconvex functions,
the sublinear function $\sigma_x(h)=\max\{p\cdot h,\; p\in\partial_-f(x)\}$.
By  a similar argument
as in the proof of \cite[Lemma 2.1.1, Chapter D]{H-UL}, we can show that $\sigma_x(h)$
satisfies \eqref{Eq.DirDer} and is therefore referred to as the directional derivative
of $f$ along $h$ \cite[Theorem 3.36]{CS}.

In the case of a locally semiconcave function $f$, we introduce the notion of
superdifferential $\partial_+f$ of $f$ at $x$ as follows.

\begin{defi}\label{Def.SupDff}
	Let $f:\Omega\mapsto \R^n$ be a locally semiconcave function in $\Omega$ and let $x\in \Omega$. Denote by $K$
	an open convex subset of $\Omega$ such that $x\in K\subset \bar{K}\subset \Omega$ and by $\omega_K$
	a semiconcave modulus for $f$ in $K$.
	The Fr\'echet superdifferential $\partial_+f$ of $f$ at $x$ is the set of vectors $p\in \R^n$ satisfying
	\begin{equation}\label{Eq.DefSupSC}
		f(y)-f(x)- p\cdot (y-x)	\leq  |y-x|\omega_K(|y-x|)
	\end{equation}
	for any point $y$ such that the segment of ends $y$ and $x$ is contained in $K$.
\end{defi}
Similar observations and properties to $\partial_-f(x)$ can be drawn for $\partial_+f(x)$.


\setcounter{equation}{0}
\section{Proofs of results}\label{Sec.Prfs}

We first prove the main results Theorem \ref{Thm.LimV} and Corollary \ref{Coro.EdgDC} by  assuming that other results hold.
Then we establish the remaining results.


\medskip
{\bf Proof of Theorem \ref{Thm.LimV}}.
\textit{Part $(i)$:}
Without loss of generality, we may assume that $x_0=0$ is a singular point and $f(0)=0$.
Let $G$ be any bounded open set such that $0\in G\subset \bar{G}\subset \Omega$ and $r>0$ be such that $\bar B_{2r}(0)\subset G$, 
and let $f$ be semiconvex in $\bar B_{2r}(0)$ with modulus $\omega_r(\cdot)$.
Given $x\in\bar B_{2r}(0)$,  $\partial_-f(x)$  is not empty, thus 
\[
	f(y)-f(x)-p_x\cdot(y-x)\geq -|y-x|\omega_r(|y-x|),\quad y,\,x\in B_{2r}(0),\; p_x\in\partial_-f(x)
\]
hence, $-f$ is upper semi-differentiable in $B_{2r}(0)$.
By the locality property (Proposition \ref{Prop.Loc}) we also have
\[
	C^u_\lambda(f_G)(x)=\lambda|x|^2-\co_{\bar B_{r}(0)}[\lambda|\cdot|^2-f](x)
\]
for $x\in \bar B_{r/2}(0)$ provided $\lambda$ is sufficiently large, and
\[
	\lim_{h\to 0}\frac{f(x)-f(0)-\sigma_0(x)}{|h|}=0,
\]
where $\sigma_0(h)=\max\{p\cdot h,\; p\in\partial_-f(0)\}$.
Note that $\partial_-f(0)$ is compact, convex and
contains more than one point since we have assumed that $0$ is a singular point.
Let $r_0>0$ be the radius of the minimal bounding sphere of $\partial_-f(0)$.
We fix $0<\epsilon <\min\{1,\, r_0\}$, then there is $0<\delta<r/2$ such that
$|f(x)-\sigma_0(x)|\leq \epsilon|x|$ whenever $x\in \bar B_\delta(0)$ as we have assumed that $f(0)=0$.
Thus for $x\in \bar B_\delta(0)$,
\[
	\sigma_0(x)-\epsilon|x|\leq f(x)\leq \sigma_0(x)+\epsilon|x|.
\]
By the locality property, we have, when $\lambda>0$ is sufficiently large,
\[
	C^u_\lambda(\sigma_0-\epsilon|\cdot|)(0)\leq C^u_\lambda(f_G)(0)\leq C^u_\lambda(\sigma_0+\epsilon|\cdot|)(0).
\]
By \eqref{Eq.Lem2.4.3}, we have
\[
	C^u_\lambda(\sigma_0+\epsilon|\cdot|)(0)\leq C^u_{(1-\epsilon)\lambda}(\sigma_0)(0)+C^u_{\epsilon\lambda}(\epsilon|\cdot|)(0)
	=\frac{r_0^2}{4(1-\epsilon)\lambda}+\frac{\epsilon}{4\lambda}
\]
hence we obtain
\[
	\lambda V_\lambda (f_G)(0)\leq \frac{r_0^2}{4(1-\epsilon)}+\frac{\epsilon}{4}.
\]
Now by \eqref{Eq.Lem2.4.1}, we have
\[
	C^u_\lambda(\sigma_0-\epsilon|\cdot|)(0)=\frac{(r_0-\epsilon)^2}{4\lambda},
\]
so that
\[
	\frac{(r_0-\epsilon)^2}{4}\leq \lambda V_\lambda(f_G)(0)\leq  \frac{r_0^2}{4(1-\epsilon)}+\frac{\epsilon}{4}.
\]
Finally we take upper and lower limits first as $\lambda\to+\infty$, then let $\epsilon\to 0+$, we obtain 
\[
	\lim_{\lambda\to +\infty}\lambda V_\lambda(f_G)(0)=r_0^2/4,
\]
which completes the proof of \textit{Part $(i)$}. \hfill\qed


\medskip
\textit{Part $(ii)$}: Let $x_0\in\Omega$ be a singular point of $f$ and let $G$ be a
bounded open convex set such that $x_0\in G\subset \bar G\subset \Omega$.
Without loss of generality, we may assume that $x_0=0$.
Since $f$ is locally semiconvex with linear modulus, we may assume that on $\bar G$,
$f(x)=g(x)-\lambda_0|x|^2$, where $g:\bar G\mapsto\R$ is convex and $\lambda_0\geq 0$ is a constant.
Clearly $\partial_-f(0)= \partial_-g(0)$. As $f(0)=g(0)$, we may further assume that $g(0)=0$. Let
$\sigma(x)=\max\{p\cdot x,\; p\in \partial_-g(0)\}$ be the sublinear function of $g$ at $0$.

\medskip
Now for every fixed $\epsilon>0$, there is a $\delta>0$ such that $|g(x)-\sigma(x)|\leq \epsilon|x|$ whenever $x\in \bar B_{\delta}(0)$. Therefore we have
\[
	\sigma(x)-\lambda_0|x|^2\leq f(x)=g(x)-\lambda_0|x|^2\leq \sigma(x)-\lambda_0|x|^2+\epsilon|x|
\]
for $x\in \bar B_{\delta}(0)$.
By the locality property, for $x\in \bar B_{\delta/2}(0)$,
and for sufficiently large $\lambda>0$, we have
\begin{equation}\label{Eq.Bnd}
	C^u_\lambda(\sigma-\lambda_0|\cdot|^2)(x)\leq C^u_\lambda(f_G)(x)\leq C^u_\lambda(\sigma+\epsilon|\cdot|-\lambda_0|\cdot|^2)(x)
\end{equation}
Now we apply Proposition \ref{Prop.UpTrC11} to $C^u_\lambda(f_G)$, then for large $\lambda>\lambda_0$,
$C^u_\lambda(f_G)\in C^{1,1}(\bar B_{\delta/2}(0))$. Let $p_\lambda=\nabla C^u_\lambda(f_G)(0)$,
we have $|p_\lambda|\leq L_G$ and $C^u_\lambda(f_G)$ is an $L_G$-Lipschitz function (see \cite[Theorem 3.12]{ZOC} and \cite[Theorem 3.5.3]{CS}) and
\[
	|C^u_\lambda(f_G)(x)-C^u_\lambda(f_G)(0)-p_\lambda\cdot x|\leq 2\lambda|x|^2
\]
for $x\in \bar B_{\delta/2}(0)$. Thus for $x\in \bar B_{\delta/2}(0)$, we have
\[
	\begin{split}
		p_\lambda\cdot x & \leq C^u_\lambda(f_G)(x)-C^u_\lambda(f_G)(0)+2\lambda|x|^2\\
				 & \leq C^u_\lambda(\sigma+\epsilon|\cdot|-\lambda_0|\cdot|^2)(x)
				       -C^u_\lambda(\sigma-\lambda_0|\cdot|^2)(0)+2\lambda|x|^2\\
				 & = I-\frac{r_0^2}{4(\lambda+\lambda_0)}+2\lambda|x|^2.
	\end{split}
\]
Here we have used the fact that
\[
	C^u_\lambda(\sigma-\lambda_0|\cdot|^2)(0)=C^u_{\lambda+\lambda_0}(\sigma)(0)
	=\frac{r_0^2}{4(\lambda+\lambda_0)}
\]
given by \eqref{Eq.Lem2.4.1} with $\epsilon=0$.
By a similar argument to that used to show \eqref{Eq.Lem2.4.3}, we also have
\[
	I=C^u_\lambda(\sigma+\epsilon|\cdot|-\lambda_0|\cdot|^2)(x)\leq
	C^u_{(1-\epsilon)\lambda}(\sigma-\lambda_0|\cdot|^2)(x)+
	C^u_{\epsilon\lambda}(\epsilon|\cdot|)(x)=J_1+J_2
\]

\noindent Now
\[
	\begin{split}
		J_1	& = C^u_{(1-\epsilon)\lambda}(\sigma-\lambda_0|\cdot|^2)(x)\\[1.5ex]
			& = C^u_{(1-\epsilon)\lambda+\lambda_0}(\sigma)(x)-\lambda_0|x|^2\\[1.5ex]
			&  =\Big(C^u_{(1-\epsilon)\lambda+\lambda_0}(\sigma)(x)-C^u_{(1-\epsilon)\lambda+\lambda_0}(\sigma)(0)+a\cdot x\Big)+
			    \Big(C^u_{(1-\epsilon)\lambda+\lambda_0}(\sigma)(0)-a\cdot x -\lambda_0|x|^2\Big)\\[1.5ex]
			&	\leq 2\Big((1-\epsilon)\lambda+\lambda_0\Big)|x|^2
				+\frac{r_0^2}{4((1-\epsilon)\lambda+\lambda_0)}-a\cdot x  - \lambda_0 |x|^2.
	\end{split}
\]
Here we have used \eqref{Eq.EstDer} and applied Lemma \ref{Lem.Ab} to the sublinear function $y\mapsto \sigma(y)$ to obtain that
 $\nabla C^u_{(1-\epsilon)\lambda+\lambda_0}(\sigma)(0)=-a$, where $-a$ is the centre of
the minimal bounding sphere of $\partial_-g(0)$, and
$C^u_{(1-\epsilon)\lambda+\lambda_0}(\sigma)(0)=r_0^2/(4((1-\epsilon)\lambda+\lambda_0))$.
We will  deal with $J_2=C^u_{\epsilon\lambda}(\epsilon|\cdot|)(x)$ later.
Therefore, when  $\lambda>\lambda_0$ is sufficiently large, we have
\[
	\begin{split}
		I-\frac{r_0^2}{4(\lambda+\lambda_0)}+2\lambda|x|^2 & \leq 2((1-\epsilon)\lambda+\lambda_0)|x|^2+\frac{r_0^2}{4((1-\epsilon)\lambda+\lambda_0)}-a\cdot x +
				C^u_{\epsilon\lambda}(\epsilon|\cdot|)(x)\\[1.5ex]
				&\phantom{xxxxxxxxxxxxxxxxx} -\frac{r_0^2}{4(\lambda+\lambda_0)}+2\lambda|x|^2-\lambda_0|x|^2\\[1.5ex]
		&	\leq \frac{\epsilon r_0^2}{4(1-\epsilon)\lambda}+8\lambda|x|^2+C^u_{\epsilon\lambda}(\epsilon|\cdot|)(x) - a \cdot x,
	\end{split}
\]
so that
\begin{equation}\label{Eq.In.01}
	(p_\lambda+a)\cdot x\leq \frac{\epsilon r_0^2}{4(1-\epsilon)\lambda}+8\lambda|x|^2+C^u_{\epsilon \lambda}(\epsilon|\cdot|)(x).
\end{equation}
Now we take
\[
	x_\lambda=\frac{p_\lambda + a}{2^5(1+|a|+L_G)\lambda},
\]
Then $|x_\lambda|\leq 1/(2^4\lambda)<\delta/2$ if $\lambda>\lambda_0$ is sufficiently large. Also
$|x_\lambda|<1/(2\lambda)$ so that
\[
	C^u_{\epsilon\lambda}(\epsilon|\cdot|)(x_\lambda)=\epsilon\lambda|x_\lambda|^2+\frac{\epsilon}{4\lambda}
\]
in the explicit formula \eqref{Eq.Lem2.4.4}. Thus if we substitute $x_\lambda$ into \eqref{Eq.In.01}, we obtain
\[
	\begin{split}
		\frac{|p_\lambda+a|^2}{2^5(1+|a|+L_G)\lambda} & \leq
				\frac{\epsilon r_0^2}{4(1-\epsilon)\lambda}+\frac{|p_\lambda+a|^2}{2^7(1+|a|+L_G)^2\lambda}+
				\frac{\epsilon\lambda|p_\lambda+a|^2}{2^{10}(1+|a|+L_G)^2\lambda^2}+\frac{\epsilon}{4\lambda}\\[1.5ex]
			& \leq \frac{|p_\lambda+a|^2}{2^7(1+|a|+L_G)\lambda}
			+\frac{\epsilon|p_\lambda+a|^2}{2^{10}(1+|a|+L_G)\lambda}+\frac{\epsilon r_0^2}{4(1-\epsilon)\lambda}+\frac{\epsilon}{4\lambda}.
	\end{split}
\]
As $0<\epsilon<1$, we have
\[
	|p_\lambda+a|^2\leq 2^6(1+|a|+L_G)(\lambda+\lambda_0)\left(\frac{\epsilon r_0^2}{4(1-\epsilon)\lambda}+\frac{\epsilon}{4\lambda}\right).
\]
Let $\lambda\to+\infty$ in the inequality above, we obtain
\[
	\limsup_{\lambda\to+\infty}|p_\lambda+a|^2\leq 2^7(1+|a|+L_G)\left(\frac{\epsilon r_0^2}{4(1-\epsilon)}+\frac{\epsilon}{4}\right).
\]
Finally, we let $\epsilon\to 0+$ and deduce that $p_\lambda\to -a$ as $\lambda\to +\infty$. Thus
\[
	\lim_{\lambda\to+\infty}\nabla C^u_\lambda(f_G)(0)=-a
\]
with $-a$ the centre of the minimal bounding sphere of $\partial_-g(0)$, which completes the proof of Part $(ii)$.
\hfill\qed.


\medskip
\begin{nota}
	We do not know whether a version of Theorem \ref{Thm.LimV}$(ii)$ holds for locally semiconvex functions with general modulus.
	To establish a similar result by following a similar approach, we need to know the regularity properties of
	$C^u_\lambda(f_G)(x)$ better in order to make the proof work.
\end{nota}


\medskip
{\bf Proof of Corollary \ref{Coro.EdgDC}:} Again, without loss of generality, we may assume that $x_0=0$ and $r_{g,0}<r_{h,0}$.
Since $E_\lambda(f_G)(0)=R_\lambda(f_G)(0)+V_\lambda(f_G)(0)\geq 0$, if $r_{g,0}=r_{h,0}$, \eqref{Eq.DCEdg} holds.
If $r_{g,0}>r_{h,0}$, as $E_\lambda(f_G)=E_\lambda(-f_G)$, we can reduce the problem to the case $r_{g,0}<r_{h,0}$.

\medskip
Next we prove, under our assumption that $r_{g,0}<r_{h,0}$ that
\begin{equation}
	\liminf_{\lambda\to\infty}\lambda R_\lambda(f_G)(0)\geq (r_{g,0}-r_{h,0})^2/4.
\end{equation}
By the locality property (see Proposition \ref{Prop.Loc}), if $\bar B_r(0)\subset G$ for some $r>0$, we see that
for $\lambda>0$ sufficiently large, we have
\[
	\co[f_G+\lambda|\cdot|^2](0)=\co_{\bar B_r(0)}[g-h+\lambda|\cdot|^2](0).
\]
Let $\sigma_g(x)=\max\{p\cdot x,\; p\in\partial_-g(0)\}$ and $\sigma_h(x)=\max\{p\cdot x,\; p\in\partial_-h(0)\}$
for $x\in\R^n$ be the sublinear functions of $g$ and $h$ at $0$ respectively,
we have, by \eqref{Eq.DirDer} that for $0<\epsilon<r_{h,0}-r_{g,0}$, there is a $0<\delta\leq r$, such that
\[
	\left|\Big(g(x)-h(x)\Big)-\Big(g(0)-h(0)\Big)-\Big(\sigma_g(x)-\sigma_h(x)\Big)\right|\leq \epsilon|x|
\]
whenever $x\in\bar B_\delta(0)$, so that
\[
	g(x)-h(x)\leq \Big(\sigma_g(x)-\sigma_h(x)\Big)+\epsilon|x|+\Big(g(0)-h(0)\Big)
\]
for $x\in\bar B_\delta(0)$. Without loss of generality, we may assume that $f(0)=g(0)-h(0)=0$.

\medskip
Again by the locality property, if $\lambda>0$ is sufficiently large, we have
\[
	\co[\lambda|\cdot|^2+f_G](0)=\co_{\bar B_\delta(0)}[\lambda|\cdot|^2+g-h](0)\leq
	\co[\lambda|\cdot|^2+\sigma_g-\sigma_h+\epsilon|\cdot|](0).
\]
Let $a_g$ be the centre of the minimal bounding sphere of $\partial_-g(0)$ and $\ell(x)=a_g\cdot x$ for $x\in\R^n$, we have
\[
\begin{split}
	\sigma_g(x)& = \max\{p\cdot x,\; p\in\partial_-g(0)\}-\ell(x)+\ell(x)\\[1.5ex]
		   & = \max\{(p-a_g)\cdot x,\;
			p\in\partial_-g(0)\}+\ell(x)\\[1.5ex]
		   & \leq r_{g,0}|x|+\ell(x)\,.
\end{split}
\]
Since the convex envelope is affine co-variant, that is $\co[H+\ell]=\co[H]+\ell$, we see that
\[
	\co[\lambda|\cdot|^2+\sigma_g-\sigma_h+\epsilon|\cdot|](0)
	\leq \co[\lambda|\cdot|^2+(r_{g,0}+\epsilon)|\cdot|-\sigma_h](0)+\ell(0).
\]
Since $\ell(0)=0$, $C^l_\lambda(H)=-C^u_\lambda(-H)$ for continuous functions $H$ of linear growth,
we may use \eqref{Eq.Lem2.4.1} in Lemma \ref{Lem.Ab} to obtain
\[
	\begin{split}
		C^l_\lambda((r_{g,0}+\epsilon)|\cdot|-\sigma_h)(0)& = \co[\lambda|\cdot|^2+(r_{g,0}+\epsilon)|\cdot|-\sigma_h](0)\\[1.5ex]
								  & = -C^u_\lambda(\sigma_h-(r_{g,0}+\epsilon)|\cdot|)(0)\\[1.5ex]
								  & = -\frac{(r_{h,0}-r_{g,0}-\epsilon)^2}{4}.
	\end{split}
\]
Thus
\[
	C^l_\lambda(f_G)(0)\leq -\frac{(r_{h,0}-r_{g,0}-\epsilon)^2}{4\lambda}
\]
when $\lambda>0$ is sufficiently large. Therefore
\[
	\lambda R_\lambda(f_G)(0)\geq \frac{(r_{h,0}-r_{g,0}-\epsilon)^2}{4}.
\]
If we let $\lambda\to+\infty$, then let $\epsilon\to 0+$, we have
\[
	\liminf_{\lambda\to+\infty}\lambda E_\lambda(f_G)(0)\geq
	\liminf_{\lambda\to+\infty}\lambda R_\lambda(f_G)(0)\geq
	\frac{(r_{h,0}-r_{g,0})^2}{4}.
\]
The proof is finished. \hfill\qed


\medskip
{\bf Proof of Proposition \ref{Prop.TghAprxScnvx}}: Suppose that $f:\Omega\mapsto\R$ is locally semiconvex with linear modulus.
Without loss of generality, we assume that $x_0=0$ is an Alexandrov point.
We set $\lambda_0=\|B\|$, the operator norm of the symmetric matrix $B$ given by \eqref{Eq.AlexTheo}.
For $\epsilon=1$, by \eqref{Eq.AlexTheo}, there is some $\delta>0$ such that
\[
	|f_G(x)-f_G(0)-p\cdot x-x^TBx|\leq \epsilon|x|^2=|x|^2
\]
whenever $x\in \bar B_{\delta}(0)$. Now we consider the affine function $\ell(x)=-f_G(0)-p\cdot x$.
Clearly $\ell(0)=-f_G(0)$. We show that
$\ell(x)\leq \lambda|x|^2-f(x)$ for all $x\in\mathbb{R}^n$ when $\lambda>0$ is large enough, so that
$-f_G(0)=\co[\lambda|\cdot|^2-f_G](0)$ hence $f_G(0)=C^u_\lambda(f_G)(0)$.

\medskip
We have, in $\bar B_{\delta}(0)$ that
\[
	-f_G(x)\geq -f_G(0)-p\cdot x-x^TBx-|x|^2\geq \ell(x)-(\lambda_0+1)|x|^2
\]
so that
\[
	\lambda|x|^2-f_G(x)\geq \ell(x)+(\lambda-\lambda_0-1)|x|^2\geq \ell(x)
\]
if $x\in \bar B_{\delta}(0)$ and $\lambda\geq \lambda_0+1$.

\medskip
If $|x|>\delta$, note that since $f_G$ is a Lipschitz function with Lipschitz constant $L_G\geq 0$, we then have
\[	
	\lambda|x|^2-f_G(x)\geq \lambda|x|^2-L_G|x|-f_G(0),
\]
while
$\ell(x)=-f_G(0)-p\cdot x\leq -f_G(0)+|p||x|$. Thus $\lambda|x|^2-f_G(x)\geq \ell(x)$ if
$\lambda|x|-L_G\geq |p|$, which holds if $\lambda\delta\geq L_G+|p|$, that is, $\lambda\geq (L_G+|p|)/\delta$.
Thus if
\[
	\lambda\geq \max\left\{\lambda_0+1,\; \frac{L_G+|p|}{\delta}\right\}
\]
we have $\lambda|x|^2-f_G(x)\geq \ell(x)$ for all $x\in \R^n$. Therefore $f_G(0)=C^u_\lambda(f_G)(0)$.

\medskip
Since in $G$, $f_G(x)=f(x)=g(x)-\lambda_1|x|^2$ for some convex function $g:\bar G\mapsto \R$ and
for some $\lambda_1>0$, if we let $\ell(x)= g(0) + q \cdot x$ for some $q\in\partial_-g(0)$, then clearly $\ell(0)=g(0)=f_G(0)$.
We show that
$g(0)+q\cdot x\leq f_G(x)+\lambda|x|^2$ for all $x\in \R^n$, hence
$f_G(0)=g(0)=\co[f_G+\lambda|\cdot|^2](0)=C^l_\lambda(f_G)(0)$
when $\lambda>0$ is sufficiently large.

\medskip
Since $0\in G$ and $G$ is open, there is a $\delta>0$ such that $\bar B_\delta(0)\subset G$. Thus in
$\bar B_\delta(0)$, we have
\[
	f_G(x)+\lambda|x|^2=g(x)+(\lambda-\lambda_1)|x|^2\geq g(x)\geq g(0)+q\cdot x
\]
if $\lambda\geq \lambda_1$.

\medskip
If $|x|>1$, similar to the proof for the upper transform, again we have $f_G(x)+\lambda|x|^2\geq \ell(x)$ 
when $\lambda>0$ is sufficiently large. Thus $f_G(0)=C^l_\lambda(f_G)(0)$ when $\lambda>0$ is sufficiently large.

\medskip
The equalities in \eqref{Eq.TghtAprx.Drv} are direct consequences of Lemma \ref{Lem.DfUp}.
Here we have $C^l_\lambda(f_G)\leq f_G\leq C^u_\lambda(f_G)$ and
$C^l_\lambda(f_G)(0)= f_G(0)\leq C^u_\lambda(f_G)(0)$, we may deduce that
$\nabla C^l_\lambda(f_G)(0)=\nabla C^u_\lambda(f_G)(0)=-p$, hence $\nabla f_G(0)=-p$.
\hfill\qed


\medskip
{\bf Proof of Lemma \ref{Lem.Ab}:}
We establish  \eqref{Eq.Lem2.4.1} first by calculating
\[
	C^u_\lambda(\sigma-\epsilon|\cdot|)(0)=-\co[\lambda|\cdot|^2+\epsilon|\cdot|-\sigma](0).
\]
We write
\[
	f_\lambda(x)=\lambda|x|^2+\epsilon|x|-\sigma(x)
\]
for $x\in\R^n$ and let $S=\partial_-f(0)$. Again let $S_{r}(-a)$ be the minimal bounding sphere of $S$
given by Lemma \ref{Lem.BndSph}. We set
\[
	b=-\frac{(r-\epsilon)^2}{4\lambda}
\]
and define the affine function $\ell(x)=a\cdot x+b$. We show that $(i)$ for $p^\ast\in S_{r}(-a)\cap S$, if we let
\begin{equation}\label{Eq.xstar} 
	x^\ast=\frac{(|p^\ast+a|-\epsilon)}{2\lambda}\frac{p^\ast+a}{|p^\ast+a|},
\end{equation}
then $f_\lambda(x^\ast)-a\cdot x^\ast=b$;
and $(ii)$ if $x^\ast$ is a minimum point of $f_\lambda(x)-a\cdot x$
then there is some $p^\ast\in S_{r}(-a)\cap S$ such that $x^\ast$ satisfies \eqref{Eq.xstar} and
$f_\lambda(x^\ast)-a\cdot x^\ast=b$.

\medskip
We prove $(i)$ first. Suppose \eqref{Eq.xstar} holds. We have
\[
	\begin{split}
		f_\lambda(x^\ast)-a\cdot x^\ast&=\lambda|x^\ast|^2+\epsilon|x^\ast|-\sigma(x^\ast)-a\cdot x^\ast\\[1.5ex]
		&=\frac{(|p^\ast+a|-\epsilon)^2}{4\lambda}-\max\{(p+a)\cdot x^\ast,\; p\in S\}+\epsilon\frac{|p^\ast+a|-\epsilon}{2\lambda}\\[1.5ex]
		&=\frac{(|p^\ast+a|-\epsilon)^2}{4\lambda}+\epsilon\frac{|p^\ast+a|-\epsilon}{2\lambda}-(p^\ast+a)\cdot x^\ast\\[1.5ex]
		&=-\frac{(|p^\ast+a|-\epsilon)^2}{4\lambda}=b.
	\end{split}
\]
Here we have used the facts that $x^\ast$ is along the direction of $p^\ast+a$ and
$p^\ast+a\in \partial (S+a)$ is the maximum point of $\max\{(p+a)\cdot x^\ast,\; p\in S\}$,
where $\partial(S+a)$ is the relative boundary of the bounded closed convex set $S+a:=\{p+a,\,p\in S\}$.

\medskip
Since $b<0$, clearly $x=0$ is not a minimum point of $f_\lambda(x)-a\cdot x$.
As the function $f_\lambda(x)-a\cdot x$ is coercive, and continuous, it reaches its minimum.
Let $x^\ast\neq 0$ be such a point.
Let $b^\prime<0$ be the minimum value of $f_\lambda(x)-a\cdot x$, that is,
$f_\lambda(x^\ast)-a\cdot x^\ast=b^\prime<0$.
Then as $-\sigma(x)$ is upper semi-differentiable and $\epsilon|x|$ is differentiable for $x\neq 0$, to follows from Lemma \ref{Lem.DfUp} that 
$\nabla (f_\lambda(x^\ast)-a\cdot x^\ast)=0$, that is
\[
	2\lambda x^\ast+\epsilon\frac{x^\ast}{|x^\ast|}-(p^\ast+a)=0
\]
where $\max\{p\cdot x^\ast,\; p\in S\}=p^\ast\cdot x^\ast$ and $p^\ast\in \partial S$,
that is, $p^\ast$ must be a relative boundary point of $S$.
Clearly, $x^\ast$ is along the same direction as $p^\ast+a$. It is easy to see that
\[
	|x^\ast|=\frac{|p^\ast+a|-\epsilon}{2\lambda}>0
\]
as $x^\ast\neq 0$. Therefore $x^\ast$ is given by \eqref{Eq.xstar}. Thus
\[
	b_\lambda=f_\lambda(x^\ast)-a\cdot x^\ast=-\frac{(|p^\ast+a|-\epsilon)^2}{4\lambda}\geq -\frac{(r_0-\epsilon)^2}{4\lambda}=b.
\]
Thus $b_\lambda=b$, hence $b=\co[f_\lambda](0)$ which implies that
\[
	\lambda V_\lambda(\sigma-\epsilon|\cdot|)(0)=
	\lambda C^u_\lambda(\sigma-\epsilon|\cdot|)(0)=-b=\frac{(r-\epsilon)^2}{4\lambda},
\]
and this proves \eqref{Eq.Lem2.4.1}.

\medskip
Now we prove \eqref{Eq.Lem2.4.2}, that is, $\nabla C^u_\lambda(\sigma)(0)=-a$ .
Let $f_\lambda(x)=\lambda|x|^2-\sigma(x)$. We have found that
$\ell(x)=a\cdot x+b\leq f_\lambda(x)$ for all $x\in\R^n$, including the special case $\epsilon=0$,
where $-a$ is the centre of the minimal bounding sphere of $\partial_-g(0)$ and $b=-r^2/(4\lambda)$.
Since $f_\lambda(x)=\lambda|x|^2-\sigma(x)$ is upper semi-differentiable in $\R^n$, by \cite{KK},
$\co[f_\lambda]\in C^1(\R^n)$. In particular $\ell(x)\leq \co[f_\lambda](x)$ and
$b=\ell(0)=\co[f_\lambda](0)$. By Lemma \ref{Lem.DfUp}, we see that $a=\nabla \ell(0)=\nabla \co[f_\lambda](0)$.
Thus by definition, $\nabla C^u_\lambda(\sigma)(0)=-a$.

\medskip
Next we establish \eqref{Eq.Lem2.4.3}. By \eqref{Eq.InLip} we have
\[
	C^u_\lambda(\sigma_0+\epsilon|\cdot|)(x)\leq C^u_{(1-\epsilon)\lambda}(\sigma_0)(x)+
	C^u_{\epsilon\lambda}(\epsilon|\cdot|)(x)
\]
for $x\in \R^n$. At $x=0$, we have, by \eqref{Eq.Lem2.4.1} with $\epsilon=0$ that
\[
	C^u_{(1-\epsilon)\lambda}(\sigma)(0)=\frac{r^2}{4(1-\epsilon)\lambda}.
\]
Also it is easy to see by a direct calculation that $C^u_{\epsilon\lambda}(\epsilon|\cdot|)(x)$ is given by \eqref{Eq.Lem2.4.4}.
Thus at $x=0$,
\[
	C^u_{\epsilon\lambda}(\epsilon|\cdot|)(0)=\frac{\epsilon}{4\lambda}\,,
\]
which completes the proof. \hfill\qed


\medskip
{\bf Proof of Proposition \ref{Prop.Loc}:}
Without loss of generality, we may assume that $x=0$. By \cite[Remark 2.1]{Z}, we have
\begin{equation}\label{Eq.LocLwTr}
	C^l_\lambda(f)(0)=\co[f+\lambda|(\cdot)-x|^2](0)=\sum^k_{i=1}\lambda_i[f(x_i)+\lambda|x_i|^2
\end{equation}
for some $2\leq k\leq n+1$, $\lambda_i>0$, $x_i\in\R^n$ for $i=1,2,\dots,k$ with $\sum^k_{i=1}\lambda_i=1$
and $\sum^k_{i=1}\lambda_ix_i=0$.  We define $f_\lambda(y)=f(y)+\lambda|y|^2$ for $y\in\R^n$.
Since $(x_i, f_\lambda(x_i))$ with $i=1,2,\dots,k$ lie on a support hyperplane of the  epi-graph
$\Epi(f_\lambda):=\{(y,\alpha),\; y\in\R^n,\; \alpha\geq f_\lambda(y)\}$, there is an affine
function $\ell(y)=a\cdot y+b$ such that
\begin{itemize}
	\item[$(i)$] $\ell(y)\leq f_{\lambda}(y)$ for all $y\in\R^n$ and
	\item[$(ii)$] $\ell(x_i)=f_\lambda(x_i)$ for $i=1,2\dots,k$.
\end{itemize}	
By $(ii)$ and \eqref{Eq.LocLwTr} we also have $\ell(0)=b=C^l_\lambda(f)(0)$. So (iii) also holds.

\medskip
To derive the bound $r_\lambda$ we evaluate $(i)$ at $y=a/(2\lambda)$ to derive a bound of $|a|$ as follows:
\[
	\frac{a\cdot a}{2\lambda}+b=\ell\left(\frac{a}{2\lambda}\right)\leq f\left(\frac{a}{2\lambda}\right)+\lambda\left|\frac{a}{2\lambda}\right|^2,
\]
so that
\[
	\frac{|a|^2}{4\lambda}\leq f\left(\frac{a}{2\lambda}\right)-b=f\left(\frac{a}{2\lambda}\right)-f(0)+f(0)-b
	\leq \frac{L|a|}{2\lambda}+\frac{L^2}{4\lambda},
\]
hence $|a|^2\leq 2L|a|+L^2$. Here we have used the fact that $f$ is $L$-Lipschitz and
$f(0)-b=R_\lambda(f)(0)\leq L^2/(4\lambda)$ by \eqref{Eq.EstLip}. Thus we have $|a|\leq (1+\sqrt{2})L$.

\medskip
Now we use $(ii)$ $a\cdot x_i+b=f(x_i)+\lambda|x_i|^2$ to obtain
\[
	\lambda|x_i|^2=b-f(x_i)+a\cdot x_i=b-f(0)+f(0)-f(x_i)+a\cdot x_i
	\leq L|x_i|+|a||x_i|,
\]
as $b-f(0)=-R_\lambda(f)(0)\leq 0$. Thus we can deduce that for each $x_i$ with $i=1,2,\dots,k$,
\[
	|x_i|\leq \frac{L+|a|}{\lambda}\leq \frac{(2+\sqrt{2})L}{\lambda}.
\]
Therefore $r_\lambda=(2+\sqrt{2})L/\lambda$. \hfill\qed


\medskip
{\bf Proof of Proposition \ref{Prop.UpTrC11}:}
We use the locality property (Proposition \ref{Prop.Loc}) to localise
the global $C^{1,1}$ property obtained in \cite[Proposition 3.7]{BKK} and
\cite[Theorem 4.1]{Z}. We show that when $\lambda>0$ is sufficiently large,
$C^u_\lambda(f)$ is continuously differentiable in $\bar B_r(0)$ and
\begin{equation}\label{Eq.DifUpTr}
	-\lambda_0|y_0-x_0|^2\leq C^u_\lambda(f)(y_0)-C^u_\lambda(f)(x_0)-\nabla
	C^u_\lambda(f)(x_0)\cdot(y_0-x_0)\leq \lambda|y_0-x_0|^2
\end{equation}
for $x_0,\; y_0\in\bar B_r(0)$, where $\lambda_0\geq 0$ is the non-negative constant used in the definition that
$f$ is semiconvex in  $\bar B_{2r}$ satisfying $f(y)=\lambda_0|y|^2-g(y)$ with $g:\bar B_{2r}\mapsto\R$ convex.
From \eqref{Eq.DifUpTr} we see that if $\lambda\geq \lambda_0$ is sufficiently large, $C^u_\lambda(f)$ is
both $2\lambda$-semiconvex and $2\lambda$-semiconcave. Therefore by \cite[Corollary 3.3.8]{CS},
$C^u_\lambda(f)\in C^{1,1}(\bar B_r(0))$ and
\[
	|\nabla C^u_\lambda(f)(y)-\nabla C^u_\lambda(f)(x)|\leq 2\lambda|y-x|
\]
for $x,\, y\in \bar B_r(0)$.

\medskip
Since $f$ is $L$-Lipschitz, by the locality property, when $\lambda>0$ is sufficiently large, we have,
for $x_0\in \bar B_r(0)$,
\[
	C^u_\lambda(f)(x_0)=-\co[\lambda|\cdot|^2-f](x_0)
	 =-\sum^{k^{(0)}}_{i=1}\lambda^{(0)}_i[\lambda|x^{(0)}_i-x_0|^2-f(x^{(0)}_i)]
\]
with $1\leq k^{(0)}\leq n+1$, $\lambda^{(0)}_i>0$, $|x^{(0)}_i-x_0|\leq r$.

\medskip
We  define $g_\lambda(y)=\lambda|y-x_0|^2-f(y)$. By Proposition \ref{Prop.Loc}, there is an affine function
$\ell(y)=a\cdot(y-x_0)+b$ such that
$(i)$: $\ell(y)\leq g_\lambda(y)$ for all $y\in\R^n$ and
$(ii)$: $\ell(x^{(0)}_i)= g_\lambda(x^{(0)}_i)$.
Let
\[
	\Delta_{x_0}=\left\{\sum^{k^{(0)}}_{i=1}\mu_i x^{(0)}_i,\;\mu_i\geq 0,\; i=1,\dots,k^{(0)},\;\sum^{k^{(0)}}_{i=1}\mu_i=1\right\}
\]
be the simplex defined by $\{x^{(0)}_1,\dots, x^{(0)}_{k^{(0)}}\}$, then we see that
$\co[g_\lambda](y)=a\cdot (y-x_0)+b$ for $y\in \Delta_{x_0}$ as the set
$U:=\{(y, a\cdot y+b),\; y\in \Delta_0\}$ is contained in a face of the convex hull of the epi-graph
$\co[\Epi(g_\lambda)]$ of $g_\lambda$
and $\{(x^{(0)}_1,g_\lambda(x^{(0)}_1)\dots,(x^{(0)}_m,g_\lambda(x^{(0)}_m)\}\subset U\cap \Epi(g_\lambda)$.

\medskip
Now we have $\co[g_\lambda](y)\leq  g_\lambda(y)$ for $y\in B_{2r}(0)$,
and $\co[g_\lambda](x^{(0)}_i)=g_\lambda(x^{(0)}_i)=a\cdot(x^{(0)}_i-x_0)+b$
for $i=1,\dots, k^{(0)}$. Furthermore, in $\bar B_{2r}(0)$,
$g_\lambda(y)=\lambda|y-x_0|^2-f(y)$ where $f(y)=g(y)-\lambda_0|y-x_0|^2$ is $2\lambda_0$-semicovex
in $\bar B_{2r}(0)$ with $g:\bar B_{2r}(0)\mapsto\R$ convex and  $\lambda_0\geq 0$.
Thus $g_\lambda(y)=(\lambda+\lambda_0)|y-x_0|^2-g(y)$ is upper semi-differentiable in $\bar B_{2r}(0)$.
Thus by Lemma \ref{Lem.DfUp}, we see that
both $\co[g_\lambda]$ and $g_\lambda$ are  differentiable at $x^{(0)}_i$ and
\[
	\nabla \co[g_\lambda](x^{(0)}_i)=\nabla g_\lambda(x^{(0)}_i)=2(\lambda+\lambda_0)(x^{(0)}_i-x_0)-\nabla g(x^{(0)}_i),
\]
hence $\nabla g(x^{(0)}_i)$ exists for $i=1,\dots, k^{(0)}$.
If we apply Lemma \ref{Lem.DfUp} to the affine function $\ell(y)$ and the upper semi-differentiable function
$g_\lambda(y)$ in $\bar B_{2r}(0)$, we also have $\nabla g_\lambda(x^{(0)}_i)=a$ for $i=1,\dots, k^{(0)}$.

\medskip
Now we show that $C^u_\lambda(f)$ is differentiable at $x_0$ and $\nabla C^u_\lambda(f)(x_0)=-a$.
We follow an argument in \cite{KK}.
We know that $\co[g_\lambda](x_0)=\sum^{k^{(0)}}_{i=1}\lambda^{(0)}_i g_\lambda(x^{(0)}_i)$ with $1\leq k^{(0)}\leq n+1$,
and we may further assume that $\lambda^{(0)}_1\geq \cdots\geq \lambda^{(0)}_{k^{(0)}}>0$,
$|x^{(0)}_i-x_0|\leq r$ (by the locality property), satisfying
$\sum^{k^{(0)}}_{i=1}\lambda^{(0)}_i=1$ and $\sum^{k^{(0)}}_{i=1}\lambda^{(0)}_ix^{(0)}_i=x_0$.
We then have $\lambda^{(0)}_1\geq 1/(n+1)$. Now for $y\in\R^n$, we have
\[
	x_0+y=\lambda^{(0)}_1\left(x^{(0)}_1+\frac{y}{\lambda^{(0)}_1}\right)+\sum^{k^{(0)}}_{i=2}\lambda^{(0)}_i x^{(0)}_i.
\]
By the convexity of $\co[g_\lambda]$, we have
\[
	\begin{split}
		\co[g_\lambda](x_0+y)-\co[g_\lambda](x_0)
		  & \leq \lambda^{(0)}_1\left(g_\lambda(x^{(0)}_1+y/\lambda^{(0)}_1)-g_\lambda(x^{(0)}_1)\right)
				+\left(\sum^{k^{(0)}}_{i=2}\lambda^{(0)}_i [g_\lambda(x^{(0)}_i)-g_\lambda (x_i^{(0)})]\right) \\[1.5ex]
		& =\lambda^{(0)}_1\left(g_\lambda(x^{(0)}_1+y/\lambda^{(0)}_1)-g_\lambda(x^{(0)}_1)\right)
	\end{split}
\]
for $y\in\R^n$. Since the left hand side of the above equation is convex in $y$ and the right hand side is
upper semi-differentiable at $y=0$ and the two terms are equal at $y=0$, by Lemma \ref{Lem.DfUp}, we see that
$\nabla \co[g_\lambda](x_0)=\nabla g_\lambda(x^{(0)}_1)$. Thus $\co[g_\lambda]$ is differentiable at $x_0$.

\medskip
Furthermore, since $\ell(y)\leq g_\lambda(y)$ for $y\in \R^n$, by the definition of convex envelope, we see that
$\ell(y)\leq \co[g_\lambda](y)$ for $y\in\R^n$. We also have $\ell(x_0)=b=\co[g_\lambda](x_0) $. Also since
$\co[g_\lambda]$ is differentiable at $x_0$, by Lemma \ref{Lem.DfUp}, we have $\nabla\co[g_\lambda](x_0)=a$.
Thus $\nabla C^u_\lambda(f)(x_0)=-a$. Therefore  $C^u_\lambda(f)$ is differentiable in $\bar B_r(0)$.
The continuity of   $\nabla C^u_\lambda(f)$ in $\bar B_r(0)$ follows from \cite{KK}.

\medskip
Now we prove that for all $x_0,\, y_0\in \bar B_r(0)$, we have
\begin{equation}\label{Eq.SubDifUpTr}
	C^u_\lambda(f)(y_0)-C^u_\lambda(f)(x_0)-\nabla C^u_\lambda(f)(x_0)\cdot (y_0-x_0)\geq -\lambda_0|y_0-x_0|^2
\end{equation}
so that  $C^u_\lambda(f)$ is $2\lambda_0$-semiconvex in $\bar B_r(0)$.
We use the notation associated to $C^u_\lambda(f)(x_0) $ as above. We see that \eqref{Eq.SubDifUpTr} is equivalent to
\[
	\lambda|y_0-x_0|^2-\co[g_\lambda](y_0)+\co[g_\lambda](x_0)+\nabla \co[g_\lambda](x_0)\cdot (y_0-x_0)\geq -\lambda_0|y_0-x_0|^2
\]
which again is equivalent to
\[
	\co[g_\lambda](y_0)-\co[g_\lambda](x_0)-\nabla \co[g_\lambda](x_0)\cdot(y_0-x_0)
	 \leq (\lambda+\lambda_0)|y_0-x_0|^2.
\]
Note that
\[
	\co[g_\lambda](x_0)
	=\sum^{k^{(0)}}_{i=1}\lambda^{(0)}_i g_\lambda(x^{(0)}_i), \quad \nabla \co[g_\lambda](x_0)=a,\quad
	\nabla \co[g_\lambda](x^{(0)}_i)=\nabla g_\lambda(x^{(0)}_i)=a.
\]
Since $y_0\in \bar B_r(0)$ and $|x^{(0)}_i-x_0|\leq r$, we see that
\[
	y_0+(x^{(0)}_i-x_0)\in \bar B_{2r}(0) \quad\text{and}\quad \sum^{k^{(0)}}_{i=1}\lambda^{(0)}_i\Big(y_0+(x^{(0)}_i-x_0)\Big)=y_0\,.
\]
Thus,
\[
	\co[g_\lambda](y_0)\leq \sum^{k^{(0)}}_{i=1}\lambda^{(0)}_i
	\co[g_\lambda](y_0+(x^{(0)}_i-x_0))\leq \sum^{k^{(0)}}_{i=1}\lambda^{(0)}_i
	g_\lambda(y_0+(x^{(0)}_i-x_0)).
\]
We also have
\[
	\co[g_\lambda](x_0)=\sum^{k^{(0)}}_{i=1} \lambda^{(0)}_i[g_\lambda(x_0+(x^{(0)}_i-x_0))]
\]
and
\[
\begin{split}
	\nabla \co[g_\lambda](x_0)\cdot(y_0-x_0)& = a\cdot (y_0-x_0)=\sum^{k^{(0)}}_{i=1}\lambda^{(0)}_ia\cdot(y_0-x_0)\\[1.5ex]
						& =\sum^{k^{(0)}}_{i=1}\lambda^{(0)}_i\nabla g_\lambda(x_0+(x^{(0)}_i-x_0)\cdot (y_0-x_0).
\end{split}
\]
We notice that in $\bar B_{2r}(0)$, $f$ is semiconvex and $f(y)=\lambda_0|y-x_0|^2-g(y)$ for a convex function
$g:\bar B_{2r}(0)\mapsto\R$. Thus
\[
	\begin{split}
		&\co[g_\lambda](y_0)-\co[g_\lambda](x_0)-\nabla \co[g_\lambda](x_0)\cdot(y_0-x_0)\\
		&\leq \sum^{k^{(0)}}_{i=1}\lambda^{(0)}_i\Big(g_\lambda(y_0+(x^{(0)}_i-x_0))-g_\lambda(x_0+(x^{(0)}_i-x_0))
			-\nabla g_\lambda(x_0+(x^{(0)}_i-x_0))\cdot(y_0-x_0)\Big)\\
		&=\sum^{k^{(0)}}_{i=1}\lambda^{(0)}_i(\lambda+\lambda_0)\Big(|(y_0-x_0)+(x^{(0)}_i-x_0))|^2-|(x^{(0)}_i-x_0))|^2
			-2(x^{(0)}_i-x_0)\cdot (y_0-x_0)\Big)\\
		&-\sum^{k^{(0)}}_{i=1}\lambda^{(0)}_i\Big(g(y_0+(x^{(0)}_i-x_0))-g(x_0+(x^{(0)}_i-x_0))
			-\nabla g(x_0+(x^{(0)}_i-x_0))\cdot (y_0-x_0)\Big)\\
		&\leq (\lambda+\lambda_0) |y_0-x_0|^2.
	\end{split}
\]
Here we have used the facts that $\sum^{k^{(0)}}_{i=1}\lambda^{(0)}_i(x^{(0)}_i-x_0)=0$ and that
$g$ is convex and differentiable at $x^{(0)}_i$. Thus $C^u_\lambda(f)$ is $2\lambda_0$-semiconvex in
$\bar B_r(0)$. Also by the definition of the upper transform, $C^u_\lambda(f)$ is $2\lambda$-semiconcave, 
hence for $x_0,\, y_0\in \bar B_r(0)$
\begin{equation}\label{Eq.LwBnSubUpTr}
	C^u_\lambda(f)(y_0)-C^u_\lambda(f)(x_0)-\nabla C^u_\lambda(f)(x_0)\cdot (y_0-x_0)\leq \lambda|y_0-x_0|^2
\end{equation}
Combining \eqref{Eq.SubDifUpTr} and \eqref{Eq.LwBnSubUpTr} we see that $C^u_\lambda(f)$ is $2\lambda_0$-semiconvex and $2\lambda$-semiconvex in $ \bar B_r(0)$.
Therefore by \cite[Corollary 3.3.8]{CS}, we see that $C^u_\lambda(f)\in C^{1,1}(\bar B_r(0))$ satisfying
\[
	|\nabla C^u_\lambda(f)(y)-\nabla C^u_\lambda(f)(x)|\leq 2\lambda|y-x|,\qquad  y,\; x\in \bar B_r(0)
\]
if we choose $\lambda\geq \lambda_0$.
\hfill\qed



\begin{thebibliography}{44}

\bibitem{A}	P. Albano,
		{\sl Some properties of semiconcave functions with general modulus},
		J. Math. Anal. Appl. {\bf 271} (2002) 217-231.
		
		\bibitem{AC}	P. Albano, P. Cannarsa,
		{\sl Structural properties of singularities of semiconcave functions},
		Ann. Scuola Norm. Sup. Pisa Cl. Sci. (4) {\bf 28} (1999) 719-740.
		
		\bibitem{AAC}	G. Alberti, L. Ambrosio, P. Cannarsa,
		{\sl On the singularities of convex functions},
		Manuscripta Math. {\bf 76} (1992) 421-435.



\bibitem{AA}	H. Attouch, D. Aze,
		{\sl Approximations and regularizations of arbitrary functions in Hilbert spaces by the Lasry-Lions methods},
		Anal. Non-Lin. H. Poincar\'e Inst. 10 (1993) 289-312.



\bibitem{BKK}	J. M. Ball, B. Kirchheim, J. Kristensen,
		{\sl Regularity of quasiconvex envelopes},
		Calc. Var. PDEs {\bf 11} (2000), 333-359.

\bibitem{Bl}	H. Blum,
		{\sl A transformation for extracting new descriptors of shape},
		Prop. Symp. Models for the Perception of Speech and Visual Form (W. W. Dunn ed.), MIT Press (1967) 362-380.


\bibitem{BW}	L. M. Blumenthal,  G. E. Wahlin,
		{\sl On the spherical surface of smallest radius enclosing a bounded subset of $n$-dimensional euclidean space},
		Bull. Amer. Math. Soc. {\bf 47}, (1941). 771-777.

\bibitem{CS}	P. Cannarsa,  C. Sinestrari,
		Semiconcave Functions, Hamilton-Jacobi Equations and Optimal Control,
		Birkh\"auser, Boston,  2004.

\bibitem{DGK}	L. Danzer, B. Gr\"unbaum, V. Klee,
		{\sl Helly's theorem and its relatives}
		In: Proceedings of Symposia in Pure Mathematics, vol. VII, pp. 101-180. AMS, Providence, RI (1963)

\bibitem{Ev}	L. C. Evans,
		Partial Differential Equations.
		Graduate Studies in Mathematics, vol. 19.
		AMS, Providence, RI, Second Edition, 2010.

\bibitem{EG}	L. C. Evans, R. F. Gariepy,
		Measure Theory and Fine Properties of Functions,
		Studies in Advanced Mathematics.
		CRC Press, Boca Raton, FL, 1992.

\bibitem{FG}	K. Fischer, B. G\"artner,
		{\sl The smallest enclosing ball of balls: combinatorial structure and algorithms},
		Int. J. Comput. Geom. Appl. {\bf 14}, (2004) 341-378.

\bibitem{Ha}	P. Hartman,
		{\sl On functions representable as a difference of convex functions},
		Pacific J. Math. {\bf 9} (1959), 707-713.

\bibitem{He}	E. Helly,
		{\sl\"Uber Mengen konvexer K\"orper mit gemeinschaftichen Punkten},
		Jber. Deutsch. Math. Verein {\bf 32}, (1923) 175–176.
		
		\bibitem{UL}	J.-B. Hiriart-Urruty,
		{\sl Generalized differentiability, duality and optimization
		for problems dealing with differences of convex functions},
		In: Convexity and duality in optimization (Groningen, 1984), 37-70,
		Lecture Notes in Econom. and Math. Systems, 256, Springer, Berlin, (1985).

\bibitem{H-UL}  J.-B. Hiriart-Urruty, C. Lemar\'echal,
		Fundamentals of Convex Analysis,
		Springer, 2001.

\bibitem{Jac92}	P. T. Jackway,
		{\sl Morphological scale-space},
		In: Proceedings 11th IAPR International Conference on Pattern Recognition.
		The Hague, The Netherlands: IEEE Computer Society Press, Los Alamitos, CA, (1992),
		pp. 252-255.

\bibitem{Jung}  H. W. E. Jung,
		{\sl \"Uber die kleinste Kugel, die eine r\"aumliche Figur einschliesst},
		J. Reine Angew. Math. {\bf 123} (1901), 241-257.

\bibitem{KK}	B. Kirchheim, J. Kristensen,
		{\sl Differentiability of convex envelopes},
		C. R. Acad. Sci. Paris Sér. I Math. {\bf 333} (2001),  725-728.

\bibitem{LL86}	J. M. Lasry and P. L. Lions,
		{\sl A remark on regularization in Hilbert Spaces},
		Israel Math. J. 55 (1986) 257-266.

\bibitem{Mor65}	J.-J. Moreau,
		{\sl Proximat\'e dualit\'e dans un espace Hilbertien},
		Bull. Soc. Math. Fr. 93 (1965) 273-299.

\bibitem{Mor66}	J.-J. Moreau,
		{\sl Fonctionnelles convexes},
		S\'{e}minaire "Sur les \'{e}quations aux d\'{e}riv\'{e}es partielles".
		Lecture Notes, Coll\'{e}ge de France, 1966.

\bibitem{OBSC}	A. Okabe, B. Boots, K. Sugihara, S. N. Chiu,
		Spatial Tessellations – Concepts and Applications of Voronoi Diagrams.
		John Wiley \& Sons, Second Edition, 2000.

\bibitem{R}	R. T. Rockafellar,
		Convex Analysis,
		Princeton Univ. Press, 1966.

\bibitem{Ser82}	J. Serra,
		Image Analysis and Mathematical Morphology.
		Volume 1, Academic Press, London, 1982


\bibitem{SP}	K. Siddiqi, S. M. Pizer (Eds),
		Medial Representations,
		Springer, New York, 2008.


\bibitem{Sy1}	J. J. Sylvester,
		{\sl A question in the geometry of situation},
		Quarterly J.  Pure and Appl. Math. {\bf 1} (1857) 79-79.

\bibitem{Sy2}	J. J. Sylvester,
		{\sl On Poncelet's approximate valuation of surd forms},
		Philosophical Magazine {\bf 20} (1860) 203-222.



\bibitem{V}	S. Verblunsky,
		{\sl On the circumradius of a bounded set},
		J. London Math. Soc. {\bf 27} (1952) 505-507.



\bibitem{Z}	K. Zhang,
		{\sl Compensated convexity and its applications},
		Anal. nonlin. H. Poincare Inst. {\bf 25} (2008) 743-771.

\bibitem{Z2}	K. Zhang,
		{\sl Convex analysis based smooth approximations of maximum functions and squared-distance functions},
		J. Nonlinear Convex Anal.  {\bf 9}  (2008) 379-406.
		
		\bibitem{ZOCc}  K. Zhang, E.C.M. Crooks, A. Orlando
		{\sl Compensated convexity, multiscale medial axis maps and sharp regularity of the squared distance function},
		SIAM J. Math. Anal., {\bf 47} (2015) 4289-4331.


\bibitem{ZOC}	K. Zhang, A. Orlando, E.C.M. Crooks,
		{\sl Compensated convexity and Hausdorff stable geometric singularity extraction},
		Math. Models Methods Appl. Sci. {\bf 25} (2015) 747-801.

 \bibitem{ZOCb}	K. Zhang, A. Orlando, E.C.M. Crooks,
		{\sl Compensated convexity and Hausdorff stable extraction of intersections for smooth manifolds},
		Math. Models Methods Appl. Sci. {\bf 25} (2015)  839-873.



\bibitem{ZOCd}  K. Zhang,  A. Orlando, E.C.M. Crooks
		{\sl UK Patent: Image Processing},
		Serial Number GB2488294, October (2015).
\end{thebibliography}
\end{document}